\documentclass[12pt]{article} \textwidth=6in \oddsidemargin=0in
\usepackage{amssymb}
\begin{document}\def\ov{\over} \def\ep{\varepsilon}
\newcommand{\C}[1]{{\cal C}_{#1}} \def\inv{^{-1}}\def\be{\begin{equation}}
\def\ee{\end{equation}}\def\x{\xi}\def\({\left(} \def\){\right)} \def\iy{\infty} \def\e{\eta} \def\cd{\cdots} \def\ph{\varphi} \def\ps{\psi} \def\La{\Lambda} \def\s{\sigma} \def\ds{\displaystyle}
\newcommand{\xs}[1]{\x_{\s(#1)}} \newcommand\sg[1]{{\rm sgn}\,#1}
\newcommand{\xii}[1]{\x_{i_{#1}}} \def\ld{\ldots} \def\Z{\mathbb Z}
\def\a{\alpha}\def\g{\gamma}\def\b{\beta} \def\noi{\noindent}
\def\ar{\leftrightarrow} \def\d{\delta} \def\S{\mathbb S} \def\tl{\tilde} \def\E{\mathbb E} \def\P{\mathbb P} \def\TS{U\backslash T}\def\CS{{\cal S}} \newcommand{\br}[2]{\left[{#1\atop #2}\right]}
\def\bs{\backslash} \def\CD{{\cal D}} \def\N{\mathbb N} \def\si{\s\inv}
\def\l{\ell}

\hfill  March 16, 2011
\begin{center}{\Large\bf Integral Formulas for the Asymmetric\\ 
\vskip1ex Simple Exclusion Process}\end{center}

\begin{center}{\large\bf Craig A.~Tracy}\\
{\it Department of Mathematics \\
University of California\\
Davis, CA 95616, USA}\end{center}

\begin{center}{\large \bf Harold Widom}\\
{\it Department of Mathematics\\
University of California\\
Santa Cruz, CA 95064, USA}\end{center}

\begin{abstract} In this paper we obtain general integral formulas for probabilities in the asymmetric simple exclusion process (ASEP) on the integer lattice $\Z$ with nearest neighbor hopping rates $p$ to the right and $q=1-p$ to the left. For the most part we consider an $N$-particle system but for certain of these formulas we can take the $N\to\iy$ limit. First we obtain, for the $N$-particle system, a formula for the probability of a configuration at time $t$, given the initial configuration. For this we use Bethe Ansatz ideas to solve the master equation, extending a result of Sch\"utz for the case $N=2$. The main results of the paper, derived from this, are integral formulas for the probability, for given initial configuration, that the $m$th left-most particle is at $x$ at time $t$. In one of these formulas we can take the $N\to\iy$ limit, and it gives the probability for an infinite system where the initial configuration is bounded on one side.  For the special case of the totally asymmetric simple exclusion process (TASEP) our formulas reduce to the known ones.
\end{abstract}

\renewcommand{\theequation}{1.\arabic{equation}}
\begin{center}{\bf I. Introduction}\end{center}

Since its introduction nearly forty years ago \cite{Spi}, the asymmetric simple exclusion
process (ASEP) has become the ``default stochastic model for transport phenomena'' \cite{Yau}.
Recall \cite{Li1, Li2} that the ASEP on the integer lattice $\Z$ is a continuous time
 Markov process $\eta_t$
where $\eta_t(x)=1$ if $x\in \Z$ is occupied at time $t$, and $\eta_t(x)=0$ if $x$ is vacant at time $t$.
Particles move on $\Z$ according to two rules: (1) A particle at $x$ waits an exponential time with
parameter one, and then chooses $y$ with probability $p(x,y)$;\linebreak
 (2) If $y$ is vacant at that time it moves
to $y$, while if $y$ is occupied it remains at $x$.  The adjective ``simple'' refers to the fact that the allowed
jumps are only one step to the right, $p(x,x+1)=p$, or one step to the left, $p(x,x-1)=q=1-p$.  The
totally asymmetric simple exclusion process (TASEP) allows jumps only to the right, so that $p=1$. 

In a major breakthrough, Johansson \cite{Jo}, building on earlier work of Baik, Deift, and Johansson
\cite{BDJ}, related a probability in TASEP to a probability in random matrix theory. Specifically, if the initial configuration in TASEP is $\Z^-$ then the probability that a particle initially at $-m$ moves at least $n$ steps to the right in time $t$ equals the probability distribution of the largest eigenvalue in a (unitary) Laguerre random matrix ensemble.  The realization \cite{Jo, PS} that the TASEP is a determinantal
process \cite{HKPV, Sos} has led to considerable progress in our understanding of the one-dimensional TASEP. (See \cite{Spo} for a recent review.)

It is natural to ask to what extent these results for TASEP can be extended to ASEP.  Since we no longer have the determinantal structure
that is present in TASEP,  random matrix theory methods, RSK-type bijections, or nonintersecting path techniques are not applicable (or at least not obviously so) to ASEP.

However, it has been known for some time \cite{ADHR, GS} that the generator of ASEP is a similarity transformation of the quantum spin chain Hamiltonian known as
the XXZ model \cite{Sut, Yang2}. Since the XXZ Hamiltonian
is diagonalizable by the Bethe Ansatz \cite{Yang2}, it is reasonable to expect that these ideas are useful for ASEP.  Indeed, Gwa and Spohn \cite{GS} applied Bethe Ansatz methods to the TASEP for a finite number of particles with periodic boundary conditions (i.e., for particles on a circle) to compute
the dynamical scaling exponent.  

Subsequently for TASEP on $\Z$ for a finite number of particles, Sch\"utz \cite{Sc}
showed that the probability that at time $t$
the system is in configuration $X=\{x_1,\ld,x_N\}$, given that its initial configuration
was $Y=\{y_1,\ldots, y_N\}$, is expressible as an $N\times N$ determinant. From this determinant representation  R\'akos and Sch\"utz \cite{RS} derived
Johansson's result relating TASEP to the Laguerre ensemble.
The R\'akos-Sch\"utz derivation uses the crucial fact that for TASEP 
the probability for any particle depends only on the initial positions for that particle and those to its right, and so it is expressible in terms of probabilities for finite systems when the initial configuration is $\Z^-$.\footnote{The determinantal formula can be use to obtain scaling limit results for other initial conditions as well. See,  e.g., \cite{Sa}.}  This is clearly no longer the case for ASEP.

In this paper we obtain general integral formulas for probabilities in ASEP. For the most part we consider an $N$-particle system but for certain of these formulas we can take the $N\to\iy$ limit, so that there are analogous formulas (involving infinite series) for infinite systems where the initial configuration is
\[Y=\{y_1,y_2,\ld\},\quad\quad y_1<y_2<\cd\to +\iy.\]
To specialize to infinite systems in TASEP with initial configuration $\Z^-$ we would replace $\Z^-$ by $Y=\Z^+$  and let $q=1$.

Denote by $P_Y(X;t)$ the probability that a system with initial configuration $Y$ is in configuration $X$ at time $t$. Then (Theorem 2.1) for an $N$-particle system  $P_Y(X;t)$ is equal to a sum of $N!$ $N$-dimensional contour integrals. The integrand was suggested by the Bethe Ansatz (but there is no Ansatz!). The case $N=2$ was established by Sch\"utz \cite{Sc} in different notation, and he also proposed that there was a general result such as this, with some contours.

The main objective of the paper is to obtain probabilities for the individual particles at time $t$. To state the formulas we introduce some notation. We set
\[\ep(\x)=p\,\x\inv+q\,\x-1,\]
and define an $N$-dimensional integrand $I(x,Y,\x)$, with variable $\x=(\x_1,\ld,\x_N)$ and parameter $x\in\Z$, by
\be I(x,Y,\x)=\prod_{i<j}{\x_j-\x_i\ov p+q\x_i\x_j-\x_i}\;{1-\x_1\cd\x_N\ov(1-\x_1)\cd(1-\x_N)}\,\prod_i\(\x_i^{x-y_i-1}e^{\ep(\x_i)t}\).\label{I}\ee
The parameter $t$ appears in the last factor but not in the notation. 
 
For the probability for $x_1(t)$, the position of the first particle at time $t$, the sum of $N!$ integrals given in Theorem~2.1 miraculously collapses into one integral. When $p\ne0$ the first product in (\ref{I}) is analytic when all the $\x_i$ lie inside some circle $\C{r}$ with center zero and radius $r$. We also assume that $r<1$. We show that for $p\ne0$ and such $r$, we have (Theorem~3.1)
\be\P(x_1(t)=x)=p^{N(N-1)/2}\int_{\C{r}}\cdots\int_{\C{r}} \,I(x,Y,\x)\,d\x_1\cd d\x_N.\label{firstPx1}\ee

In order to take the $N\to\iy$ limit we need integrals over large contours rather than small ones. This is because of the factors $\x_i^{-y_i}$ in the integrand (\ref{I}). Given a set $S\subset\{1,\ld,N\}$ we denote by $I(x,Y_S,\x)$ the integrand analogous to $I(x,Y,\x)$ where only the variables $\x_i$ with $i\in S$ occur. Then 
when $q\ne0$ we have (Theorem~3.2)
\be\P(x_1(t)=x)=\sum_S c_S\,\int_{\C{R}}\cd\int_{\C{R}} I(x,Y_S,\x)\,d^{|S|}\x,\label{Px1}\ee
where $R>1$ is so large that all the poles of the first product in the integrand lie inside $\C{R}$. The sum runs over all nonempty subsets $S$ of $\{1,\ld,N\}$ and $c_S$ are certain constants involving $S$ and powers of $p$ and $q$.

Once we have this we are able to compute the expected value $\E(x_1(t))$ and to take the $N\to\iy$ limit in (\ref{Px1}). If 
$Y=\{y_1,\,y_2,\ld\}$ with
$y_1<y_2<\cd\to+\iy$ 
then on the right side of (\ref{Px1}) we simply take the sum over all finite sets $S\subset\Z^+$, the resulting series being convergent.

These results are contained in Section~III. In Section~IV we derive the analogue of (\ref{firstPx1}) for the second-left particle.

The main results of the paper are in Section V where we obtain the analogues of (\ref{firstPx1}) and (\ref{Px1}) for the general particle. The analogue of (\ref{firstPx1}) has the form
\be\P(x_{m}(t)=x)=\sum_{|S|>N-m}c_{S,m,N}\int_{\C{r}}\cd\int_{\C{r}}I(x, Y_S,\x)\,d^{|S|}\x,\label{Pxm1}\ee
where we sum over all subsets $S$ of $\{1,\ld,N\}$ with cardinality at least $N-m+1$, and $c_{S,m,N}$ is another explicitly given constant. (See Theorem~5.1.) The analogue of (\ref{Px1}) is
\be\P(x_{m}(t)=x)=\sum_{|S|\ge m}c_{S,m}\int_{\C{R}}\cd\int_{\C{R}}I(x, Y_S,\x)\,d^{|S|}\x,\label{Pxm2}\ee
where we sum over all subsets $S$ of $\{1,\ld,N\}$ with cardinality at least $m$, and $c_{S,m}$ is yet another explicitly given constant. (See Theorem~5.2.) Notice that in the latter representation the coefficients are independent of $N$, and this allows us to take the $N\to\iy$ limit and so obtain probabilities for infinite systems. When the initial configuration is $\Z^+$ the sum over sets in (\ref{Pxm2}) may be simplified to a sum over~integers. (The corollary to Theorem 5.2.)

Of course (\ref{firstPx1}) and (\ref{Px1}) and the second-particle formula are special cases of these but we give their proofs first, so that  we can introduce the new ingredients gently. 

The deduction of these formulas from Theorem~2.1 requires two algebraic identities which were discovered by computer computation of special cases. The first is
\[\sum_{\s\in{\S}_N} {\rm sgn}\,\s\,\left(\prod_{i<j}\,(p+q\x_{\s(i)}\x_{\s(j)}-\x_{\s(i)})\right.\]
\[\left.\times{\xs{2}\,\x_{\s(3)}^2\,\x_{\s(4)}^3\cdots \x_{\s(N)}^{N-1}\ov
(1-\x_{\s(1)}\x_{\s(2)}\x_{\s(3)}\cdots\x_{\s(N)})\cdots (1-\x_{\s(N-1)}\x_{\s(N)})(1-\x_{\s(N)})}\right)\]
\be=p^{N(N-1)/2} {\prod_{i<j} (\x_j-\x_i)\ov \prod_j (1-\x_j)},\label{identity1}\ee
where the sum is over all permutations $\s$ in the symmetric group $\S_N$.
We also use an equivalent version of this identity,
\pagebreak
\[
\sum_{\s\in\S_N} \sg\s\left(\prod_{i<j}\left( p+ q\x_{\s(i)}\x_{\s(j)}-\x_{\s(i)}\right)\right.\]
\[\left. \times{1\ov (\x_{\s(1)}-1)(\x_{\s(1)}\x_{\s(2)}-1)\cdots (\x_{\s(1)}\x_{\s(2)}\cdots \x_{\s(N)}-1)}\right)\]
\be=  q^{N(N-1)/2}\,{\prod_{i<j}(\x_j-\x_i)\ov \prod_j(\x_j-1)},\label{identity1'} \ee
obtained from (\ref{identity1}) by interchanging $p$ and $q$ and letting $\x_i\to 1/\x_{N-i+1}$.

For the second identity we introduce the notation
\[[N]={p^N-q^N\ov p-q}\,,\]
and then
\be[N]!=[N]\,[N-1]\,\cdots\,[1]\,,\quad \br{N}{m}={[N]!\ov [m]! \;[N-m]!}\,,\label{qbinom}\ee
where we set $[0]!=1$. (Note that $\br{N}{m}$ is $q^{m\,(N-m)}$ times a $q$-binomial coefficient with $q$ equal to our $p/q$. Hence the notation.) The identity is
\be\sum_{|S|=m} \prod_{{ i\in S\atop  j\in S^c}} {p+q\x_i\x_j-\x_i\ov \x_j-
\x_i}\cdot\Big(1-\prod_{j\in S^c}\x_j\Big)=q^m\,\br{N-1}{m}\, \Big(1-\prod_{j=1}^N \x_j\Big)\label{identity2}\ee
for $N\ge m+1$. The sum runs over all subsets $S$ of $\{1,\ld,N\}$ with cardinality $m$, and $S^c$ denotes the complement of $S$ in $\{1,\ld,N\}$. 

The proofs of these identities will be given in the last section.

\setcounter{equation}{0}\renewcommand{\theequation}{2.\arabic{equation}}
\begin{center}{\bf II. Solution of the Master Equation}\end{center}

We denote by $Y=\{y_1,\ld,y_N\}$ with $y_1<\cd<y_N$ the initial configuration of the process and write $X=\{x_1,\ld,x_N\}\in\Z^N$. When $x_1<\cd<x_N$ then $X$ represents a possible configuration of the system at a later time $t$. We denote by $P_Y(X;t)$ the probability that the system is in configuration $X$ at time $t$, given that it was initially in configuration $Y$.

Given $X=\{x_1,\ld,x_N\}\in\Z^N$ we set
\[X_i^+=\{x_1,\ld,x_{i-1},x_i+1,x_{i+1},\ld,x_N\},\quad
X_i^-=\{x_1,\ld,x_{i-1},x_i-1,x_{i+1},\ld,x_N\}.\]
The master equation for a function $u$ on $\Z^N\times\mathbb{R}^+$ is
\be{d\ov dt}u(X;t)=\sum_{i=1}^N\Big(p\,u(X_i^-;t)+q\,u(X_i^+;t)-u(X;t)\Big),\label{master}\ee
and the boundary conditions are, for $i=1,\ld,N-1$,
\pagebreak
\[u(x_1,\ldots,x_i,x_i+1,\ldots,x_N;t)\]
\be = p\, u(x_1,\ld, x_i,x_i,\ld,x_N;t)+q\, u(x_1,\ld, x_i+1,x_i+1,\ldots,x_N;t).\label{boundary}\ee
The initial condition is
\[u(X;0)=\d_Y(X)\ \ {\rm when}\ x_1<\cd<x_N.\]
The basic fact is that if $u(X;t)$ satisfies the master equation, the boundary conditions, and the initial condition, then $P_Y(X;t)=u(X;t)$
when $x_1<\cd<x_N$.\footnote{The idea in Bethe Ansatz (see, e.g., \cite{LL, Sut, Yang2}), applied to one-dimensional
$N$-particle quantum mechanical problems, is to represent the wave function as a linear combination of free particle eigenstates and to incorporate the effect of the potential as a
set of $N-1$ boundary conditions.  The remarkable
feature of models amendable to Bethe Ansatz is that the boundary conditions for $N\ge 3$ introduce
no more new conditions, with the result that (\ref{boundary}) involves only consecutive particles. The application of Bethe Ansatz to the evolution equation (master equation) describing ASEP begins with  Gwa and Spohn \cite{GS} with subsequent developments by  Sch\"utz \cite{Sc}.} 

Recall that an inversion in a permutation $\s$ is an ordered pair $\{\s(i),\s(j)\}$ in which $i<j$ and $\s(i)>\s(j)$. We define
\[S_{\a\b}=-{p+q\x_\a\x_\b-\x_\a\ov p+q\x_\a\x_\b-\x_\b}\]
and then
\[A_\s=\prod\{S_{\a\b}:\{\a,\b\}\ {\rm is\ an\ inversion\ in}\ \s\}.\]
We also set
\[\ep(\x)=p\,\x\inv+q\,\x-1.\] For now we shall assume $p\ne 0$, so the $A_\s$ are analytic at zero in all the variables. Here and later all differentials $d\x$ incorporate the factor $(2\pi i)\inv$. 

\noi{\bf Theorem 2.1}. We have when $p\ne0$
\be P_Y(X;t)=\sum_{\s\in\S_N}\int_{\C{r}}\cdots\int_{\C{r}} A_\s\,\prod_i\x_{\s(i)}^{x_i-y_{\s(i)}-1}\,e^{\,\sum_i\ep(\x_i)\,t}\,d\x_1\cd d\x_N,\label{theorem}\ee
where $\C{r}$ is a circle centered at zero with radius $r$ so small that all the poles of the $A_\s$ lie outside $\C{r}$.
 
\noi{\bf Remark}. For TASEP with $p=1$ we have 
\[S_{\a\b}=-(1-\x_\a)/(1-\x_\b).\]
 Fix $\a$. In $A_\s$ the factor $1-\x_\a$ occurs for each inversion of the form $\{\a,\b\}$ and $(1-\x_\a)\inv$ for each inversion of the form $\{\b,\a\}$. If $\a=\s(i)$ then the number of the former minus the number of the latter equals $\s(i)-i$. The number of inversions is the number of transpositions whose product is $\s$. These give
\[A_\s={\rm sgn}\,\s\prod(1-\x_{\s(i)})^{\s(i)-i}.\] 
Now the integrand in (\ref{theorem}) factors and we obtain the sum of $\sg\s$ times
\[\prod_i\int_{\C{r}} (1-\xs{i})^{\s(i)-i}\,\xs{i}^{x_i-y_{\s(i)}-1}e^{\ep(\xs{i})t}\,d\xs{i}=\prod_i\int _{\C{r}}(1-\x)^{\s(i)-i}\,\x^{x_i-y_{\s(i)}-1}e^{\ep(\x)t}\,d\x,\]
and so the sum over $\s$ equals
\[\det\Big(\int_{\C{r}}(1-\x)^{j-i}\,\x^{x_i-y_{j}-1}e^{\ep(\x)t}\,d\x\Big).\]
This is the determinant representation of $P_Y(X;t)$ obtained in \cite{Sc}.  
 
To prove the theorem we shall show three things:

\noi(a) The right side of (\ref{theorem}) satisfies the master equation for all $X\in\Z^N$.

\noi(b) The right side of (\ref{theorem}) satisfies the boundary conditions for all $X\in\Z^N$.

\noi(c) The right side of (\ref{theorem}) satisfies the initial condition when $x_1<\cd<x_N$.

\noi{\bf Proof of (a)}. This is clear once the last factor in (\ref{theorem}) is written as the exponential of $\sum_i\ep(\xs{i})\,t$.

\noi{\bf Proof of (b)}. We shall show that the boundary condition is satisfied pointwise by the integrand. Let $T_i\s$ denote $\s$ with the entries $\s(i)$ and $\s(i+1)$ interchanged. The boundary conditions will be satisfied provided that\footnote{The Bethe Ansatz proposes an integrand as in (\ref{theorem}) with some coefficients $A_\s$. The fact that (\ref{boundary}) only involves consecutive particles implies that satisfying it only requires relations between the coefficients $A_s$ and $A_{T_i\s}$. A straightforward computation shows that the relations are as stated.} 
\[A_{T_i\s}=S_{\s(i+1),\s(i)}\,A_\s\]
for all $\s$.
Let us see why this relation holds. Let $\a=\s(i),\ \b=\s(i+1)$, and suppose $\a>\b$. Then $\{\a,\b\}$ is an inversion for $\s$ but not for $T_i\s$, so $S_{\a\b}$ is a factor in $A_\s$ but not in $A_{T_i\s}$, and all other factors are the same. Therefore, using $S_{\a\b}\,S_{\b\a}=1$, we have
\[A_{T_i\s}=S_{\b\a}\,A_{\s}=S_{\s(i+1),\s(i)}\,A_\s.\]
The same identity holds immediately if $\b>\a$, since $\{\b,\a\}$ is an inversion for $T_i\s$ but not for $\s$. Thus, (b) is established.

\noi{\bf Proof of (c)}. The initial condition is satisfied by the summand in (\ref{theorem}) coming from the identity permutation $id$. So what we have to show is
\[\sum_{\s\ne{\rm id}}\int_{\C{r}}\cdots\int_{\C{r}} A_\s\,\prod_i\xs{i}^{x_i-y_{\s(i)}-1}\,d\x_1\cd d\x_N=0,\]
or equivalently
\[\sum_{\s\ne{\rm id}}\int_{\C{r}}\cdots\int_{\C{r}} A_\s\,\prod_i\x_i^{x_{\si(i)}-y_{i}-1}\,d\x_1\cd d\x_N=0,\]
when $x_1<\cd<x_N$. We write $I(\s)$ for the integral corresponding to $\s$, so that the above becomes
\be\sum_{\s\ne{\rm id}}I(\s)=0.\label{Isum}\ee

For $1\le n<N$ fix $n-1$ distinct numbers $i_1,\ld,i_{n-1}\in[1,\,N-1]$. Define 
\[A=\{i_1,\ld,i_{n-1}\},\]
and then
\[\S_N(A)=\{\s\in\S_N:\s(1)=i_1,\ld,\s(n-1)=i_{n-1},\,\s(n)=N\}.\]
When $n = 1$ this consists of all permutations with $N$ in position 1, and when $n=N-1$ it consists of a single permutation. If $B$ is the complement of 
$A\cup\{N\}$  in $[1,\,N]$,
then $\s\in\S_N(A)$ is determined by the restriction
\[\s|_{[n+1,\,N]}:[n+1,\,N]\to B.\] 

\noi{\bf Lemma 2.1}. For each $A$,
\[\sum_{\s\in\S_N(A)}I(\s)=0.\]
\noi{\bf Start of the Proof}. When $\s\in\S_N(A)$ the inversions involving $N$ are the $(N,i)$ with $i\in B$. Therefore the integrands in $I(\s)$ involving these $\s$ may be written
\[\prod_{i\in B}S(\x_N,\x_i)\times
\prod_{i\le N}\x_i^{x_{\si(i)}-y_i-1}
\times\prod\{S(\x_k,\x_\l):N>k>\l,\ \s\inv(k)<\s\inv(\l)\}.\] 
The integrals are taken over $\C{r}$ with $r$ so small that all the denominators in the $S$-factors are nonzero on and inside the contour. In these integrals we make the substitution
\[\x_N\to{\e\ov\prod_{i<N}\x_i},\]
so that $\e$ runs over a circle of radius $r^N$. The integrand becomes
\be(-1)^{N-n}\prod_{i\in B}{p+q\e\prod_{\l\ne i,\,N}\x_\l\inv-
\e\,\prod_{\l\ne N}\x_\l\inv\ov p+q\e\prod_{\l\ne i,\,N}\x_\l\inv-
\x_i}\label{first}\ee
\be\times\ \e^{x_n-y_N-1}\,\prod_{i<N}\x_i^{x_{\s\inv(i)}-x_n+y_N-y_i-1}\label{second}\ee
\be \times\ \prod\{S(\x_k,\x_\l):N>k>\l,\ \s\inv(k)<\s\inv(\l)\}.\label{third}\ee
(The reason that we still have $-1$ in the exponents in (\ref{second})
is that $d\x_N=\prod_{i<N}\x_i\inv\,d\e$.)

\noi{\bf Sublemma 2.1}. When $n=N-1$ we have $I(\s)=0$.

\noi{\bf Proof}. There is a single $i\in B$ and (\ref{first})  is analytic inside the $\x_i$-contour except for a simple pole at $\x_i=0$. The power of $\x_i$ in (\ref{second}) is 
\[\xi_i^{x_N-x_{N-1}+y_N-y_i-1},\]
and since $x_N>x_{N-1}$ and $y_N>y_{i}$, the exponent is positive. Therefore the integrand is analytic inside the $\x_i$-contour, and so the integral is zero.   

\noi{\bf Sublemma 2.2}. When $n>N-1$ all $I(\s)$ with $\s\in\S_N(A)$ are sums of lower-order integrals in each of which (\ref{first}) is replaced by a factor depending on $A$. The other factors remain the same. In each integral some $\x_i$ with $i\in B$ is equal to another $\x_j$ with $j\in B$.

\noi{\bf Proof}. We may assume that $q\ne0$. This case follows by a limiting argument. We are going to shrink some of the $\x_i$-contours with $i\in B$. Due to the defining property of $r$, the only poles we pass will come from the product (\ref{first}). In fact, to avoid double poles later we take $\x_i\in \C{r_i}$ with the $r_i$ all slightly different.

Take $j=\max B$ and shrink the $\x_j$-contour. The product 
(\ref{first}) has a simple pole at $\x_j=0$ (the $j$-factor has the pole and the $i$-factors with $i\ne j$ are analytic there) and the power of $\x_j$ in (\ref{second}) is positive as before, so the integrand is analytic at $\x_j=0$. For each $k\in B$ with $k\ne j$ we pass the pole at
\be\x_j={q\e\prod_{\l\ne j,k,N}\x_\l\inv\ov \x_k-p}.\label{polej}\ee
coming from the $k$-factor in (\ref{first}). (Our assumption on the $r_i$ assures that there are no double poles.) For the residue we replace the $k$-factor by
\be -{p+q\e\prod_{\l\ne k,\,N}\x_\l\inv-
\e\,\prod_{\l\ne N}\x_\l\inv\ov q\e\x_j^{-2}\prod_{\l\ne j,k,N}\x_\l\inv},\label{k}\ee
where in this and the $j$-factor we replace $\x_j$ by the right side of (\ref{polej}). When $i\ne j,k$ the $i$-factor becomes
\[p+q\e\prod_{\l\ne i,\,N}\x_\l\inv-
\e\,\prod_{\l\ne N}\x_\l\inv\ov p\,(1-\x_i\x_k\inv),\]
and we replace $\x_j$ in the numerator by  the right side of (\ref{polej}).

We now shrink the $\x_k$-contour. There is a pole of order 2 at 
$\x_k=0$ coming from (\ref{k}) and the $j$-factor in (\ref{first}). Since $k<j=\max B<N$, we have $y_N-y_k\ge 2$, so the exponent of 
$\x_k$ in (\ref{second}) is at least 2. Therefore the integrand is analytic at $\x_k=0$. The factor (\ref{k}) has no other poles inside $\C{r_k}$. An $i$-factor with $i\ne j,k$ will have a pole at $\x_k=\x_i$  if $r_i<r_k$.  There is also the pole at
\[\x_k={q\e\prod_{\l\ne j,k,N}\x_\l\inv\ov \x_j-p}\]
coming from the $j$-factor. But this relation and (\ref{polej}) imply $\x_j=\x_k$. 

Thus when we shrink the $\x_j$-contour and the $\x_k$-contours with $k\ne j$ we obtain $(N-2)$-dimensional integrals in each of which two of the $\x$-variables corresponding to indices in $B$ are equal. This proves the sublemma. 

\noi{\bf Sublemma 2.3}. For each integral of Sublemma 2.2 there is a partition of $\S_N(A)$ into pairs $\s,\,\s'$ such that $I(\s)+
I(\s')=0$ for each pair.

\noi{\bf Proof}. Consider an integral in which $\x_i=\x_j$. We pair $\s$ and $\s'$ if $\s\inv(i)={\s'}\inv(j)$ and $\s\inv(j)={\s'}\inv(i)$, and $\s\inv(k)={\s'}\inv(k)$ when $k\ne i,j$. The factor (\ref{second}) is clearly the same for both when $\x_i=\x_j$, and we shall show that the $\s$- and $\s'$-factors in (\ref{third}) are negatives of each other when $\x_i=\x_j$.

Assume for definiteness that
\be i<j\ \ {\rm and}\ \ \s\inv(i)<\s\inv(j).\label{assume}\ee
(Otherwise we reverse the roles of $\s$ and $\s'$.) Then the factor 
$S(\x_j,\x_i)$ does not appear for $\s$ in (\ref{third}) but it does appear for $\s'$.  
This factor equals $-1$ when $\x_i=\x_j$. 

To complete the proof it is enough to show that for any $k\ne i,j$ the product of $S$-factors involving $k$ and either $i$ or $j$ is the same for $\s$ and $\s'$ when $\x_i=\x_j$. There are nine cases, depending on the position of $k$ relative to $i$ and $j$ and the position of 
$\s\inv(k)$ relative to $\s\inv(i)$ and $\s\inv(j)$. 
If $k$ is outside the interval $[i,\,j]$ and $\s\inv(k)$ is outside the interval $[\s\inv(i),\,\s\inv(j)]$ then the products of $S$-factors for $\s$ and $\s'$ are clearly the same. There are five remaining cases, with the results displayed in the table below. The first column gives the position of $k$ relative to $i$ and $j$, the second column gives the position of $\s\inv(k)$ relative to 
$\s\inv(i)$ and $\s\inv(j)$, the third column gives the product of $S$-factors involving $k$ and either $i$ or $j$ for $\s$, and the fourth column gives the corresponding product for $\s'$. Keep (\ref{assume}) in mind.

\[\begin{array}{llllll}
i<k<j&\si(k)<\si(i) && S(\x_k,\x_i) && S(\x_k,\x_j)\\
i<k<j&\si(i)<\si(k)<\si(j) && 1 && S(\x_k,\x_j)\,S(\x_i,\x_k)\\
i<k<j&\si(j)<\si(k) && S(\x_j,\x_k) && S(\x_i,\x_k)\\
k<i&\si(i)<\si(k)<\si(j) && S(\x_i,\x_k) && S(\x_j,\x_k)\\ 
k>j&\si(i)<\si(k)<\si(j) && S(\x_k,\x_j) && S(\x_k,\x_i)
\end{array}\]

In all cases but the second the $S$-factors are exactly the same for $\s$ and $\s'$ when $\x_i=\x_j$. For the second we use $S(\x_k,\x)\,S(\x,\x_k)=1$. 

Sublemmas 2.1--2.3 give Lemma 2.1. 

To prove (\ref{Isum}) we use induction on $N$. When $N=2$ it folows from Sublemma 2.1. Assume $N>2$ and that the result holds for $N-1$. For those permutations for which $\s(N)=N$ we use the fact that $N$ appears in no involution, so we may integrate with respect to $\x_1,\ld,\x_{N-1}$ and use the induction hypothesis. The set of permutations with $\s(N)<N$ is the disjoint union of the various $\S_N(A)$, and for these we apply Lemma 2.1.

This completes the proof of (c), and so of Theorem~2.1.

\noi{\bf Remark}. In case $q\ne0$ the same formula (\ref{theorem}) holds when the circle is sufficiently large instead of small. A similar argument to the one just given should hold, but there is another way to see this. If we set 
\[X^-=\{-x_N,\ld,-x_1\},\quad Y^-=\{-y_N,\ld,-y_1\}\]
and denote by $\tl P$ the probability density for the process with $p$ and $q$ interchanged, then $P_Y(X;t)=\tl P_{Y^-}(X^-;t)$. If we apply (\ref{theorem}) to this other process and then make the substitutions $\x_i\to \x_i\inv$ in the integrals we obtain (\ref{theorem}) with a large $\C{R}$. This duality will be use again in Section~V, in the derivation of (\ref{Pxm2}) from (\ref{Pxm1}).

\setcounter{equation}{0}\renewcommand{\theequation}{3.\arabic{equation}}
\begin{center}{\bf III. The Left-most Particle}\end{center}

Here we determine the probability that the left-most particle $x_1$ is at $x$ at time $t$. 

\noi{\bf Theorem 3.1}. With $\C{r}$ as before and $I(x,Y,\x)$ given by (\ref{I}), we have when $p\ne0$
\be\P(x_1(t)=x)=p^{N(N-1)/2}\int_{\C{r}}\cdots\int_{\C{r}} \,I(x,Y,\x)\,d\x_1\cd d\x_N.\label{Th2}\ee

\noi{\bf Proof}. Since $x_1<\cd<x_N$ we may rewrite $X=\{x_1,\,x_2,\ld,x_N\}$ as
\[x,\;x+z_1,\;x+z_1+z_2,\ld,\;x+z_1+\cd+z_{N-1}.\]
Then $\P(x_1(t)=x)$ equals the sum of $P_Y(X;t)$ over all $z_i>0$. After summing, the integrand in (\ref{theorem}) becomes 
\[A_\s\;{(1-\x_{\s(1)}\cd\x_{\s(N)})\,
\x_{\s(2)}\x_{\s(3)}^2\cd\x_{\s(N)}^{N-1}\ov (1-\x_{\s(1)}\x_{\s(2)}\cd\x_{\s(N)}) (1-\x_{\s(2)}\cd\x_{\s(N)})\cd(1-\x_{\s(N)})}\;\prod_i\(\x_i^{x-y_i-1}e^{\ep(\x_i)t}\).\]
If we observe that
\be A_\s={\rm sgn}\,\s\prod_{{i<j\atop\s(i)>\s(j)}}{p+q\x_{\s(i)}\x_{\s(j)}-\x_{\s(i)}
\ov p+q\x_{\s(i)}\x_{\s(j)}-\x_{\s(j)}} =\sg\s\,{\ds{\prod_{i<j}}(p+q\x_{\s(i)}\x_{\s(j)}-\x_{\s(i)})
\ov \ds{\prod_{i<j}}(p+q\x_i\x_j-\x_i)},\label{reps}\ee
then we see that the theorem follows from (\ref{identity1}).

We derive here the alternative expression for $\P(x_1(t)=x)$. 
Given a set $S\subset\{1,\ld,N\}$ there is a corresponding set $Y_S=\{y_i:i\in S\}$ and corresponding $|S|$-dimensional integrand $I(x,Y_S,\x)$. We define $\s(S)=\sum_{i\in S}i$, the sum of the indices in $S$.

\noi{\bf Theorem 3.2}. When $q\ne0$ we have 
\be \P(x_1(t)=x)=\sum_S {p^{\,\s(S)-|S|}\ov q^{\,\s(S)-|S|(|S|+1)/2}}\,\int_{\C{R}}\cd\int_{\C{R}} I(x,Y_S,\x)\,d^{|S|}\x,\label{Th3}\ee
where $R$ is so large that all the poles of the integrand lie inside $\C{R}$. The sum runs over all nonempty subsets $S$ of $\{1,\ld,N\}$. 

We begin with a lemma that replaces integrals such as appear in (\ref{Th2}) by sums of integrals over large contours. 

Suppose $f(\x_1,\ld,\x_N)$ is analytic for all $\x_i\ne0$ and that for $i>k$
\[f(\x_1,\ld,\x_N)\Big|_{\x_i\to (\x_k-p)/q\x_k}=O(\x_k)\]
as $\x_k\to0$, uniformly when all $\x_j$ with $j\ne k$ are bounded and bounded away from zero. Define
\[I_{f}(\x)=\prod_{i<j}{\x_j-\x_i\ov p+q\x_i\x_j-\x_i}\;{f(\x_1,\ld,\x_N)\ov \prod_{i}(1-\x_i)}.\]
For a subset $S$ of $\{1,\ld,N\}$ let $I_{f,S}(\x)$ denote the analogous function where the variables are the $\x_i$ with $i\in S$, and in $f(\x_1,\ld,\x_N)$ the $\x_i$ with $i\in S^c$ are replaced by~1.

\noi{\bf Lemma 3.1}. Under the  stated assumptions on $f(\x_1,\ld,\x_N)$, we have when $p,\,q\ne0$  
\be \int_{\C{r}}\cd\int_{\C{r}} I_{f}(\x)\,d^N\x=\sum_S {p^{\,|S^c|-\s(S^c)}\ov q^{\,\s(S)-|S|(|S|+1)/2}}\,\int_{\C{R}}\cd\int_{\C{R}} I_{f,S}(\x)\,d^{|S|}\x\,,\label{intsum}\ee
where $r$ is so small that the poles of the integrand on the left lie outside $\C{r}$ and $R$ is so large that the poles of the integrand on the right lie inside $\C{R}$. The sum runs over all subsets $S$ of $\{1,\ld,N\}$. When $S$ is empty the integral is interpreted as $f(1,\ld,1)$.

\noi{\bf Proof}. We use induction. The result is easily seen to be true when $N=1$, so we assume $N>1$ and that the lemma holds for $N-1$. We expand the $\x_N$-contour on the left side. In addition to the pole at $\x_N=1$ we encounter poles at $\x_N=(\x_k-p)/q\x_k$. We claim that the residue at this pole, when integrated over $\x_k$, will give zero. The factor $f(\x_1,\ld,\x_N)$, after substituting for $\x_N$ its value at the pole, is $O(\x_k)$ as $\x_k\to0$ by the assumption on $f$, while
\[{\prod_{i<j}(\x_i-\x_j)\ov \prod_i(1-\x_i)},\]
after substituting for $\x_N$ its value at the pole, will be of the order $\x_k^{-N+2}$ at $\x_k=0$. The residue of $1/(p+q\x_k\x_N-\x_k)$ at the pole is $1/q\,\x_k$. The factor $1/(p+q\x_i\x_N-\x_i)$ with $i\ne k$ equals $\x_k/(p(\x_k-\x_i))$. The factor $\x_k-\x_i$ is cancelled by the same factor in the numerator, so no new poles in the $\x_k$ variable are introduced. So the $1/(p+q\x_i\x_N-\x_i)$ combined, including the contribution of the residue, give the power $\x_k^{N-3}$. Thus the product of all factors combined is $O(1)$. Hence the $\x_k$ integral equals zero, as claimed.

So after we expand the $\x_N$-contour we have an integral where $\x_N$ is over an arbitrarily large contour $\C{R}$ and the other $\x_i$  over small contours $\C{r}$, and another integral (coming from the pole at $\x_N=1$) in which $\x_N$ does not appear and the other $\x_i$ are over $\C{r}$. 

Let us consider the latter. The integral we get is the left side of (\ref{intsum}) with $N$ replaced by $N-1$ times
\be\prod_{i<N}{\x_N-\x_i\ov p+q\x_i\x_N-\x_i}\Big|_{\x_N=1}={1\ov p^{N-1}}.\label{factor}\ee
(Notice that $\x_{N}$ appears in the denominator of $I_f(\x)$ as $1-\x_{N}$, and the pole at $\x_{N}=1$ was outside the contour. The two minus signs cancel when we compute the contribution from the pole.) Our induction hypothesis tell us that this equals
\[{1\ov p^{N-1}}\;\sum_{S\subset\{1,\ld,N-1\}} {p^{\,|S^c|-\s(S^c)}\ov q^{\,\s(S)-|S|(|S|+1)/2}}\,\int_{\C{R}}\cd\int_{\C{R}} I_{f,S}(\x)\,d^{|S|}\x.\]
Now $S^c$ here indicates complement with respect to $\{1,\ld,N-1\}$. If we want to express this in terms of the complement in $\{1,\ld,N\}$ as in the statement of the lemma we have to make the substitutions
\[|S^c|\to|S^c|-1,\quad \s(S^c)\to \s(S^c)-N,\]
and therefore in the notation of the lemma the above equals
\[\sum_{S\subset\{1,\ld,N-1\}} {p^{\,|S^c|-\s(S^c)}\ov q^{\,\s(S)-|S|(|S|+1)/2}}\,\int_{\C{R}}\cd\int_{\C{R}} I_{f,S}(\x)\,d^{|S|}\x.\]
This is the portion of the right side of (\ref{intsum}) corresponding to those $S$ not containing~$N$.

Now for the original integral, but where $\x_N$ is taken over $\C{R}$. If we integrate first with respect to $\x_N$, leaving us with an $N-1$-dimensional integral with small $\C{r}$, we see that we are in the case $N-1$ with $f(\x_1,\ld,\x_N)$ replaced by
\be\tl{f}(\x_1,\ld,\x_{N-1})=\int_{\C{R}} \prod_{i<N}{\x_N-\x_i\ov p+q\x_i\x_N-\x_i}\;{f(\x_1,\ld,\x_N)\ov 1-\x_N}\,d\x_N.\label{newf}\ee
We have to show that $\tl{f}$ satisfies the required condition. For notational convenience we take $k=1$ and $i=2$, so we make the substitution $\x_2\to(\x_1-p)/q\x_1$. All but the product are $O(\x_1)$ as $\x_1\to0$, by assumption on $f$.  Write the product as the product over $i\ne 1,2$, which is bounded if $R$ was chosen large enough, uniformly for the $\x_i$ bounded away from zero, times
\[{\x_N-\x_1\ov p+q\x_N\x_1-\x_1}\,{\x_N-(\x_1-p)/q\x_1\ov p+\x_N\,(\x_1-p)/\x_1-(\x_1-p)/q\x_1}\]
\[={\x_N-(\x_1-p)/q\x_1\ov p+q\x_N\x_1-\x_1}\,{\x_N-\x_1\ov p+\x_N\,(\x_1-p)/\x_1-(\x_1-p)/q\x_1}.\]
The first factor equals $1/q\x_1$ while the second factor is $O(\x_1)$. Thus (\ref{newf}) satisfies the required condition and we may again apply the induction hypothesis. 

If $\tl{S}\subset\{1,\ld,N-1\}$ to compute $I_{\tl f,\tl{S}}$ we replace the $\x_i$ with $i\in\tl{S}^c$ by 1 in $\tl{f}(\x_1,\ld,\x_{N-1})$. We see that the product in the integrand in (\ref{newf}) is replaced by
\[{1\ov q^{|\tl{S}^c|}}\,\prod_{i\in\tl{S}}{\x_N-\x_i\ov p+q\x_i\x_N-\x_i}.\]
If we set $S=\tl{S}\cup\{N\}$ then in terms of $S$ the full integrand including the $\x_N$-variable is 
\[{1\ov q^{|{S}^c|}}\,I_{f,S}(\x).\]
For the coefficient on the right side of (\ref{intsum}) with $S$ replaced by $\tl{S}$, the power of $p$ is unchanged while the power of $q$ is
\[\s(\tl{S})-|\tl{S}|\,(|\tl{S}|+1)/2=\s(S)-N-(|S|-1)\,|S|/2,\]
and if we add to this $|{S}^c|=N-|S|$ we get
\[\s(S)-|S|\,(|S|+1)/2,\]
which is the power of $q$ in (\ref{intsum}). 

Summing over these $S$ gives the portion of the right side of (\ref{intsum}) corresponding to those $S$ containing $N$, and this completes the proof of Lemma~3.1.

\noi{\bf Proof of Theorem 3.2}. First assume $p\ne0$. Observe that  $I(x,Y,\x)$ is $I_{f}(\x)$ as defined above with
\be f(\x_1,\ld,\x_N)=\Big(1-\prod_i\x_i\Big)\,\prod_i\(\x_i^{x-y_i-1}\,e^{\ep(\x_i)\,t}\),\label{f}\ee
and $I(x,Y_S,\x)$ is $I_{f,S}(\x)$. We show that $f$ as defined this way satisfies the hypothesis of the lemma.
The exponential summands $\ep(\x_i)$ and $\ep(\x_k)$ combine when $\x_i=(\x_k-p)/q\x_k$ to give
\[{pq\x_k\ov\x_k-p}+q{\x_k-p\ov q\x_k}+{p\ov\x_k}+q\x_k-2=q\x_k+{pq\x_k\ov\x_k-p}-1,\]
which is analytic at $\x_k=0$. The powers of $\x_i$ and $\x_k$ combine as
\[\({\x_k-p\ov q\x_k}\)^{x-y_i-1}\,\x_k^{x-y_k-1}=O(\x_k^{y_i-y_k})\]
as $\x_k\to0$, which is $O(\x_k)$ since $y_i>y_k$. So the hypothesis on $f$ is satisfied. Since $f(1,\ld,1)=0$ the sum may be taken over nonempty subsets $S$. Because of the factor $p^{N(N-1)/2}$ in (\ref{Th2}) we must multiply the factor in (\ref{intsum}) by this, resulting in the factor in (\ref{Th3}). 

We can remove the requirement $p\ne0$ by taking the $p\to0$ limit since the power of $p$ is nonnegative and no pole tends to infinity as $p\to0$.

An immediate conclusion from Theorem 3.2 is that $\P(x_1(t)=x)$ tends exponentially to zero as $x\to-\iy$, and therefore so does $\P(x_1(t)<x)$. Thus $\P(x_1(t)\ge x)$ tends exponentially to 1. 

We write
\[Q(x)=p^{N(N-1)/2}\int_{\C{r}}\cd\int_{\C{r}}\x_1^{x-y_1-1}\cd\x_N^{x-y_N-1}\,\prod_{i<j}{\x_i-\x_j\ov p+q\x_i\x_j-\x_i}\]
\[\times{1\ov\prod_i(1-\x_i)}\,e^{\,\sum_i\ep(\x_i)\,t}\,d\x_1\cd d\x_N,\]
where $r$ is small. Clearly $Q(x)\to 0$ as $x\to+\iy$ and
\[\P(x_1(t)=x)=Q(x)-Q(x+1).\]
 It follows that 
$Q(x)=\P(x_1(t)\ge x)$, and this tends exponentially to 1 as $x\to-\iy$.
Therefore
\[\E(x_1(t))=\lim_{x\to-\iy}\sum_{y=x}^\iy y\,[Q(y)-Q(y+1)]=\lim_{x\to-\iy}\left[\sum_{y=x}^\iy Q(y)+(x-1)\,Q(x)\right]\]
\[=\lim_{x\to-\iy}\left[\sum_{y=x}^\iy Q(y)+(x-1)\right].\]
Now 
\[\sum_{y=x}^\iy Q(y)=p^{N(N-1)/2}\int_{\C{r}}\cd\int_{\C{r}}\x_1^{x-y_1-1}\cd\x_N^{x-y_N-1}\,\prod_{i<j}{\x_j-\x_i\ov p+q\x_i\x_j-\x_i}\]
\be\times{1\ov 1-\x_1\cd\x_N}\cdot{1\ov\prod_i(1-\x_i)}\,e^{\,\sum_i\ep(\x_i)\,t}\,d\x_1\cd d\x_N.\label{Qint}\ee

If we apply the procedure of the last section to this integral we start by moving the $\x_N$-contour out. The resulting integral, with $\x_N$ over a very large contour, is exponentially small as $x\to-\iy$, even though the other contours are small. As before the residues at the poles $\x_N=(\x_k-p)/q\x_k$ integrate out to zero. The contribution from the pole $\x_N=1$ gives the same expression we started with but with $N$ replaced by $N-1$. But there is now also a pole at $\x_N=1/\x_1\cd\x_{N-1}$ whose contribution is 
\[-\int_{\C{r}}\cd\int_{\C{r}}\ps_{N-1}(\x_1,\ldots,\x_{N-1};y_N-y_1,\ldots,y_N-y_{N-1})\,d\x_1\cd d\x_{N-1},\]
where
\[\ps_{N}(\x_1,\ldots,\x_{N};z_1,\ldots,z_{N})=p^{N(N+1)/2}\x_1^{z_1}\cd\x_{N}^{z_N}\,\prod_{i<j}{\x_j-\x_i\ov p+q\x_i\x_j-\x_i}\]
\[\times\prod_{i}{(\x_1\cd\x_{N})\inv-\x_i\ov p+q\x_i(\x_1\cd\x_{N})\inv-\x_i}\,\cdot\,{1\ov 1-\x_1\cd\x_{N}}\,\cdot\,{1\ov\prod_{i}(1-\x_i)}\,e^{\,\left[\sum_{i}\ep(\x_i)+\ep((\x_1\cd\x_{N})\inv)\right]\,t}.\]

For the right side of (\ref{Qint}) when $N=1$, which is what we are left with at the end, we expand the contour, encounter a double pole at $\x_1=1$, and find that it equals
\[(p-q)\,t+y_1-x+1+{\rm exponentially\ small\ term}.\]
Therefore 
\[\E(x_1(t))=(p-q)\,t+y_1-\sum_{j=1}^{N-1}\int_{\C{r}}\cd\int_{\C{r}}\ps_{j}(\x_1,\ldots,\x_j;y_{j+1}-y_1,\ldots,y_{j+1}-y_{j})\,d\x_1\cd 
d\x_{j}.\]

The integral $\int_{\C{r}}\ps_1(\x;z)\,d\x$ has an explicit expression in terms of Bessel functions $I_n(2t)$. Indeed, it equals
\[2pt \,e^{-2t}\,[I_{z-1}(2t)+I_{z}(2t)]+(2z-1)p\,\Big\{{1\ov 2}\, e^{-2t} I_0(2t)-{1\ov 2}+e^{-2t}\sum_{j=1}^{z-1} I_j(2t)\Big\}.\]
The integrals of the other $\ps_j(\x;z)$ are not so simple.

We now show how to obtain the probability $\P(x_1(t)=x)$ for a system with infinitely many particles\footnote{It follows from the fact that ASEP is a Feller process \cite{Li1} that the limit as $N\to\iy$ equals the probability for the infinite system.  We thank Thomas Liggett \cite{Li3} for explaining this fact to us.} with initial configuration
\be Y=\{y_1,y_2,\ld\},\quad\quad y_1<y_2<\cd\to+\iy\label{initial}\ee
when $q\ne0$. In fact, it is very easy once we have Theorem 3.2. We just modify the right side so that the sum runs over all finite subsets of $\Z^+$. Because $\C{R}$ may be taken arbitrarily large the resulting series converges, as we shall now show.

Consider the various factors in the integrand in (\ref{Th2}), where the $\x_i$ are replaced by the $\x_n$ with $n\in S$. The $-1$ parts of the exponents of the $\x_n$ may be removed if we replace $d\x_n$ by $d\x_n/\x_n$, which is $O(1)$. Suppose that $|S|=k$ as above. The analogue of the factor
\[{1-\x_1\cd\x_N\ov \prod_i(1-\x_i)}\]
is at most $(R^k+1)/(R-1)^k=O(2^k)$ for $R$ large. The product $\prod_{i<j}(\x_i-\x_j)$ is at most $(2R)^{k(k-1)/2}$. The denominator
$\prod_{i<j}(p+q\x_i\x_j-\x_i)$ is at least $(qR^2/2)^{k(k-1)/2}$. The product of the $\x_n$ is at most $R^{kx-\sum y_n}$. The exponential factor is 
$O(e^{aRk})$ for some $a$ depending on $t$. So the integral is 
\[O\((a/R)^{k(k-1)/2}\,R^{kx-\sum y_n}\,e^{aRk}\)=R^{O(k)-\sum y_n}\]
if we take $R>a$. (Since $R$ is fixed the factor $e^{aRk}$ can be incorporated into the $R^{O(k)}$ term.) Since $\sum y_n\ge\sum (y_1+n-1)=\s(S)+k\,(y_1-1)$ the above is at most $R^{O(k)-\s(S)}$. And since $\s(S)\ge k(k+1)/2$, this is at most $R^{-\s(S)+O(\s(S)^{1/2})}$.

The external factor in (\ref{Th3}) is
\[p^{\s(S)-k}\,q^{k(k+1)/2-\s(S)}\le q^{-\s(S)}.\]
It follows that if we take $R>1/q^2$ then the integral times the external factor is at most $R^{-\s(S)/2}$. 

Now consider all sets $S$ with $\s(S)=k$. Since the largest $i\in S$ is at most $k$, the number of such sets is at most $2^{k}$. Hence the sum of the absolute values of the terms of the infinite series is at most a constant times
\[\sum_{k=1}^\iy 2^k\,R^{-k/2},\]
which is finite when $R>4$. Thus we have shown convergence for all $t$.

\setcounter{equation}{0}\renewcommand{\theequation}{4.\arabic{equation}}

\begin{center}{\bf IV. The Second-Left Particle}\end{center}

In this section we compute the probability $\P(x_2(t)=x)$. It is somewhat more complicated than that given for $\P(x_1(t)=x)$ in Theorem 3.1 and the proof introduces some new elements. 

We use the notation (\ref{I}), and for $1\le k\le N$ we set $Y_k=Y\backslash{y_k}$.

\noi{\bf Theorem 4.1}. With the contours $\C{r}$ as in Theorem 3.1 we have when $p\ne0$
\[ \P(x_2(t)=x)=-q\,{p^{N-1}-q^{N-1}\ov p-q}\,p^{(N-1)(N-2)/2}\,\int_{\C{r}}\cd\int_{\C{r}} I(x,Y,\x)\,d^N\x\]
\be +\,p^{(N-1)(N-2)/2}\sum_{k=1}^N \({q\ov p}\)^{k-1}\int_{\C{r}}\cd\int_{\C{r}} I(x,Y_k,\x)\,d^{N-1}\x.\label{Th4}\ee

\noi{\bf Proof}. We rewrite $X=\{x_1,\,x_2,\ld,x_N\}$ as
\[x-v,\;x,\;x+z_1,\ld,\;x+z_1+\cd+z_{N-2}.\]
Then $\P(x_2(t)=x)$ equals the sum of $P_Y(X;t)$ over all $v>0$ and $z_i>0$. If we sum first over  $z_2,\cd,z_{N-2}$ the product $\prod_i\x_{\s(i)}^{x_i-y_{\s(i)}-1}$ in (\ref{theorem}) becomes, 
$\prod_i\x_i^{x-y_i-1}$ times
\[\x_{\s(1)}^{-v}{\x_{\s(3)}\,\x_{\s(4)}^2\cd\x_{\s(N)}^{N-2}\ov
(1-\x_{\s(3)}\x_{\s(4)}\cd\x_{\s(N)})\cd(1-\x_{\s(N)}\x_{\s(N-1)})(1-\x_{\s(N)})}.\]

We now move the $\x_{\s(1)}$-contour out beyond the unit circle, and we do not encounter any poles. Here is the reason. From the first part of (\ref{reps}) we see that we get poles at
\be\x_{\s(i)}=\left\{\begin{array}{ll}(\x_{\s(k)}-p)/ q\x_{\s(k)}&{\rm if}\ k>i,\\&\\p/(1-q\x_{\s(k)})&{\rm if}\ k<i.\end{array}\right.\label{poles2}\ee
The poles in the $\x_{\s(1)}$-variable are at $\x_{\s(1)}=(\x_{\s(k)}-p)/q\x_{\s(k)}$, and these are very large since $\xs{k}\in\C{r}$. So we move the $\x_{\s(1)}$-contour out beyond the unit circle and sum, giving
\be{1\ov \x_{\s(1)}-1}{\x_{\s(3)}\,\x_{\s(4)}^2\cd\x_{\s(N)}^{N-2}\ov
(1-\x_{\s(3)}\x_{\s(4)}\cd\x_{\s(N)})\cd(1-\x_{\s(N)}\x_{\s(N-1)})(1-\x_{\s(N)})}.\label{term1}\ee
If we move the contour back to $\C{r}$ we pass a pole at $\xs{1}=1$ with residue
\be{\x_{\s(3)}\,\x_{\s(4)}^2\cd\x_{\s(N)}^{N-2}\ov
(1-\x_{\s(3)}\x_{\s(4)}\cd\x_{\s(N)})\cd(1-\x_{\s(N)}\x_{\s(N-1)})(1-\x_{\s(N)})}.\label{residue}\ee
The factor $A_\s$ when $\x_{\s(1)}=1$ equals
\be\sg\s\,\(-{q\ov p}\)^{\s(1)-1}\prod\Big\{{p+q\x_{\s(i)}\x_{\s(j)}-\x_{\s(i)}\ov p+q\x_{\s(i)}\x_{\s(j)}-\x_{\s(j)}}:i<j,\ \s(i)>\s(j),\ i,j>1\Big\},\label{term2}\ee
since there are $\s(1)-1$ indices less than $\s(1)$.

Now we add all terms in which $\s(1)=k$. For the contribution from the pole at $\xs{1}=1$, when $\s(1)=k$ (\ref{term2}) may be written (as in (\ref{reps}))
\[\sg\s\,\(-{q\ov p}\)^{k-1}\,{1\ov \ds{\prod_{i<j}}(p+q\x_i\x_j-\x_i)}\,\prod_{i<j}(p+q\x_{\s(i)}\x_{\s(j)}-\x_{\s(i)}),\]
where all indices are $\ne k$. This is to be multiplied by (\ref{residue}). If we use identity (\ref{identity1}) to sum the product over those $\s$ with $\s(1)=k$ we get
\[\({q\ov p}\)^{k-1}p^{(N-1)(N-2)/2}\(1-\ds{\prod_{j\ne k}}\x_j\)\,{\ds{\prod_{i<j}}(\x_j-\x_i)\ov\ds{\prod_j{(1-\x_j)}}\,{\ds{\prod_{i<j}}}(p+q\x_i\x_j-\x_i)},\]
where again all indices are $\ne k$. The exterior factor is now $\prod_{i\ne k}\(\x_{i}^{x-y_i-1}\,e^{\ep(\x_i)t}\)$, and from these we obtain Äthe sum on the right side of (\ref{Th4}).

Next we consider (\ref{term1}), which we rewrite as
\[{(1-\x_{\s(2)}\x_{\s(3)}\cd\x_{\s(N)})\ov\x_k-1}{\x_{\s(3)}\,\x_{\s(4)}^2\cd\x_{\s(N)}^{N-2}\ov
(1-\x_{\s(2)}\x_{\s(3)}\cd\x_{\s(N)})\cd(1-\x_{\s(N)}\x_{\s(N-1)})(1-\x_{\s(N)})},\]
where $\s$ is now a map from $\{2,\ldots,N\}$ to $\{1,\ldots,N\}\backslash\{k\}$. The factor ${\rm sgn}\,\s$ becomes $(-1)^{k+1}\,{\rm sgn}\,\s$. We also rewrite $A_\s$ as
\[{1\ov \ds{\prod_{i<j}}(p+q\x_i\x_j-\x_i)}\,\prod_{j\ne k}(p+q\x_k\x_j-\x_k)\,\prod_{i<j}(p+q\x_{\s(i)}\x_{\s(j)}-\x_{\s(i)}).\]
If we use identity (\ref{identity1}) to sum this over these $\s$ we get
\be(-1)^{k}\,p^{(N-1)(N-2)/2}\,\(1-\ds{\prod_{j\ne k}}\x_j\)\,{\ds{\prod_{{i<j\atop i,j\ne k}}}(\x_j-\x_i)\ov\ds{\prod_j{(1-\x_j)}}}\,{\ds{\prod_{j\ne k}}(p+q\x_k\x_j-\x_k)\ov \ds{\prod_{i<j}}(p+q\x_i\x_j-\x_i)}
.\label{Nsummand}\ee
The other factors are still $\prod_i\(\x_{i}^{x-y_i-1}\,e^{\,\sum_{i}\ep(\x_i)t}\)$.

To evaluate the sum of (\ref{Nsummand}) over $k$ we write it as
\[ -p^{(N-1)(N-2)/2}\,{\ds{\prod_{i<j}(\x_j-\x_i)}\ov \ds{\prod_j (1-\x_j)\,\prod_{i<j}(p+q\x_i\x_j-\x_i)}}\]
\[\times
\sum_k  \(1-\prod_{j\ne k}\x_j \)\,{\ds{\prod_{j\ne k}(p+q\x_j\x_k-\x_k)}\ov \ds{\prod_{j\ne k}(\x_j-\x_k)}}.\]
The factor $(-1)^k$ became $-1$ because of the way we rewrote the product of the $\x_j-\x_i$ in (\ref{Nsummand}). Identity (\ref{identity2}) with $k=1$ tells us that the last sum equals 
\[q\,{p^{N-1}-q^{N-1}\ov p-q}\,\(1-\prod_{j}\x_j \).\]
If we recall the power of $p$ in the first factor above and the ubiquitous factor $\prod_{i}\(\x_{i}^{x-y_i-1}\,e^{\ep(\x_i)t}\)$ we see that this gives the first term in (\ref{Th4}).

\setcounter{equation}{0}\renewcommand{\theequation}{5.\arabic{equation}}
 
\begin{center}{\bf V. The General Particle}\end{center}

In this section we consider the $m$th particle from the left for general $m$. We prove, with the notation (\ref{qbinom}),

\noi{\bf Theorem 5.1}. We have when $p\ne0$
\[\P(x_{m}(t)=x)=p^{(N-m)(N-m+1)/2}\,q^{m(m-1)/2}\]
\be\times\sum_{|S^c|<m}(-1)^{m-1-|S^c|}\br{|S|-1}{m-|S^c|-1}\,{q^{\s(S^c)- m\,|S^c|}\ov p^{
\s(S^c)-|S^c|(|S^c|+1)/2}}\int_{\C{r}}\cd\int_{\C{r}}I(x, Y_S,\x)\,d^{|S|}\x,\label{Th5.1}\ee
where $r$ is so small that the poles of the integrand lie outside $\C{r}$. The sum is taken over all subsets $S$ of $\{1,\ld,N\}$ with $|S^c|<m$.

We shall first establish a preliminary form of the result, and for that we use a lemma analogous to Lemma~3.1 and which follows from it.
Now we have a function $g(\x_1,\ld,\x_N)$ which is analytic for all $\x_i\ne0$ and satisfies, for $i<k$,
\[g(\x_1,\ld,\x_N)\Big|_{\x_i\to p/(1-q\x_k)}=O(\x_k\inv)\]
as $\x_k\to\iy$, uniformly when all $\x_j$ with $j\ne k$ are bounded and bounded away from zero. Define, as before,
\[I_g(\x)=\prod_{i<j}{\x_j-\x_i\ov p+q\x_i\x_j-\x_i}\;{g(\x_1,\ld,\x_N)\ov \prod_{i}(1-\x_i)}.\]
and similarly $I_{g,S}(\x)$.

\noi{\bf Lemma 5.1}. Under the stated assumptions on $g(\x_1,\ld,\x_N)$, we have when $p,\,q\ne0$
\be \int_{\C{R}}\cd\int_{\C{R}} I_g(\x)\,d^N\x=\sum_S (-1)^{|S^c|}{q^{\,\s(S^c)-N\,|S^c|}\ov p^{\,|S|(N+1+|S^c|)/2-\s(S)}}\,\int_{\C{r}}\cd\int_{\C{r}} I_{g,S}(\x)\,d^{|S|}\x\,,\label{intsum2}\ee
where $r$ is so small that the poles of the integrand on the right lie outside $\C{r}$ and $R$ is so large that the poles of the integrand on the left lie inside $\C{R}$. As before, $S$ runs over all subsets of $\{1,\ld,N\}$.

\noi{\bf Proof}. We apply Lemma 3.1 with $p$ and $q$ interchanged to the function
\[f(\x_1,\ld,\x_N)=g(\x_N\inv,\ld,\x_1\inv)\prod\x_i\inv.\]
The required hypothesis on $f$ follows from the hypothesis on $g$.  The left side of (\ref{intsum2}), after the substitutions $\x_i\to1/\x_{N-i+1}$, equals the left side of (\ref{intsum}) after interchanging $p$ and $q$, times $(-1)^N$, because $\prod(1-\x_i)$  becomes $\prod(\x_i-1)$. So we apply Lemma 3.1 to $f$ and then change variables again resulting in a factor $(-1)^{|S|}$ in each summand. The reason for the different coefficients is the reversal of the order of the variables. Thus the coefficient of the integral involving $I_{f,S}$ with $p$ and $q$ interchanged equals the coefficient of the integral involving $I_{g,\tl S}$, where $\tl S=\{N-i+1:i\in S\}$. This, together with the resulting power of $-1$, is what we see on the right side of (\ref{intsum2}). This proves Lemma~5.1.

To state the preliminary form of the result we introduce more notation.
For disjoint subsets $T$ and $U$ of $\{1,\ld,N\}$, we define
\pagebreak
 
\[I(x, Y_{T,U},\x)\]
\be =\Big(1-\prod_{i\in U} \x_i\Big)
{\ds{\prod_{{i<j\atop i,j\in U\ {\rm or}\ i,j\in T}}(\x_j-\x_i)} 
\ov\ds{\prod_{i}(1-\x_i)}}\;
{\ds{\prod_{i\in T,\;j\in U}(p+q\x_i\x_j-\x_i)}
\ov \ds{\prod_{i<j}(p+q\x_i\x_j-\x_i)}}\,\prod_{i} \(\x_i^{x-y_i-1}\,e^{\ep(\x_i)\,t}\),\label{YTU}\ee
where indices with unspecified range run over $T\cup U$.
If $T\subset U$ we define $\s(T,U)$ to be the sum of the positions of the elements of $T$ in $U$. Thus, if $U=\{2,\,3,\,5\}$ and $T=\{2,\,5\}$ then $\s(T,U)=1+3=4$. Finally, for a set $U$ we define
\[\sg U=(-1)^{\#\{(i,j)\,:\,i>j,\;i\in U,\,j\in U^c\}}.\]

\noi{\bf Lemma 5.2}. With $r$ small enough and $p,q\ne0$ we have
\[\P(x_{m}(t)=x)=p^{(N-m)(N-m+1)/2+m(m-1)/2}\,q^{(m-1)(m-2)/2}\]\[\times\sum_{|U|=m-1}\sg U\,\sum_{T\subset U}(-1)^{|U\bs T|+\s(\TS)-\s(\TS,\,T)}{q^{\s(\TS)-(m-1)\,|\TS|}\ov p^{\s(\TS)+|T|(m+|\TS|)/2}}\]
\[\times \int_{\C{r}}\cd\int_{\C{r}}I(x, Y_{T,U^c},\x)\,d^{|T\cup U^c|}\x,\]
where $U$ runs over all subsets of $\{1,\ld,N\}$ with $|U|=m-1$.
 
\noi{\bf Proof}. To begin we now write $X$ as
\[x-v_{m-1}-\cd-v_1,\;x-v_{m-2}-\cd-v_1,\;\cd, x-v_1,\;x,\;x+z_1,\ld,\;x+z_1+\cd+z_{N-m}.\]
The product $\prod_i\(\x_{\s(i)}^{x_i-y_{\s(i)}-1}\,e^{\ep(\x_i)\,t}\)$ in (\ref{theorem}) is replaced by $\prod_i\(\x_i^{x-y_i-1}e^{\ep(\x_i)t}\)$ times
\[\x_{\s(1)}^{-v_1-\cd-v_{m-1}}\cd\x_{\s(m-1)}^{-v_1}\;\x_{\s(m+1)}^{z_1}\cd\x_{\s(N)}^{z_1+\cd+z_{N-m}},\]
and we have to sum over all $v_i>0,\ z_i>0$. As in the proof of Theorem~4.1 we can sum over the $z_i$ immediately and we get
\[\x_{\s(1)}^{-v_1-\cd-v_{m-1}}\cd\x_{\s(m-1)}^{-v_1}\,
{\x_{\s(m+1)}\x_{\s(m+2)}^2\cd\x_{\s(N)}^{N-m}\ov (1-\x_{\s(m+1)}\cd\x_{\s(N)})\cd(1-\x_{\s(N)})}.\]

Before we can sum over the $v_i$ we have to move the $\x_{\s(i)}$-contours out, and we do them in the order $i=1,\ldots,m-1$. As in the proof of Theorem~4.1 we see by referring to (\ref{poles2}) that the poles obtained from moving the $\x_{\s(1)}$-contour are very large and so we can move that contour out almost that far. Then if we want to move the $\x_{\s(2)}$-contour out we encounter poles with $k>2$ in (\ref{poles2}), which are far out and so no problem, but also the pole with $k=1$ when $\s(2)<\s(1)$, which is at 
$\x_{\s(2)}=p/(1-q\x_{\s(1)})$.
We show that the residue at this pole, when integrated with respect to 
$\x_{\s(1)}$, gives zero. 

Recall that the $\x_{\s(1)}$-contour is large, and so $\x_{\s(2)}$ as a function of $\x_{\s(1)}$ is analytic outside the $\x_{\s(1)}$-contour and is in fact $O(\x_{\s(1)}\inv)$ at infinity. The part of the residue that comes from the factor 
\[{p+q\x_{\s(1)}\x_{\s(2)}-\x_{\s(1)}
\ov p+q\x_{\s(1)}\x_{\s(2)}-\xs{2}}\]
is $O(1)$ at infinity as are all the other factors in the first part of (\ref{reps}) because $\x_{\s(2)}$, in terms of $\xs{1}$, is small when $\x_{\s(1)}$ is large. The powers of the $\x_i$ involving $\x_{\s(1)}$ and $\x_{\s(2)}$ combine as
\[\x_{\s(1)}^{-v_1-\cd-v_{m-1}-x-1}\;\({p\ov1-q\x_{\s(1)}}\)^{-v_1-\cd-v_{m-2}-x-1}\,\xs{1}^{-y_{\s(1)}}\,\({p\ov1-q\x_{\s(1)}}\)^{-y_{\s(2)}},\]
which is analytic outside the $\x_{\s(1)}$-contour and $O(\x_{\s(1)}^{-v_{m-1}-y_{\s(1)}+y_{\s(2)}})$ at infinity. The exponent is $\le-2$ since $v_{m-1}>0$ and $y_{\s(2)}<y_{\s(1)}$.  Finally we have to check the  exponential of $\sum\ep(\x_i)$. The sum of those involving $\x_{\s(1)}$ and $\x_{\s(2)}$ is
\[{p\ov\xs{1}}+q\xs{1}+ (1-q\xs{1})+q{p\ov1-q\xs{1}},\]
which is bounded at infinity. Hence the $\x_{\s(1)}$-integral is zero, as claimed.

Continuing this way we move all the $\x_{\s(i)}$-contours out for $i<m$. We then sum over the $v_i$, obtaining 
\[{1\ov (\xs{1}-1)\,(\xs{1}\xs{2}-1)\cd(\xs{1}\xs{2}\cd\xs{m-1}-1)}\]
\[\times{\x_{\s(m+1)}\x_{\s(m+2)}^2\cd\x_{\s(N)}^{N-m}\ov (1-\x_{\s(m+1)}\cd\x_{\s(N)})\cd(1-\x_{\s(N)})}\]
as replacement for what we had before. We also have, from the numerator in the second part of (\ref{reps}), a product which we write as
\[\prod_{i<j<m}(p+q\x_{\s(i)}\x_{\s(j)}-\x_{\s(i)})\;\prod_{m\le i<j}(p+q\x_{\s(i)}\x_{\s(j)}-\x_{\s(i)})\;
\prod_{{i<m\le j}}(p+q\x_{\s(i)}\x_{\s(j)}-\x_{\s(i)}).\]
Therefore we are to take the sum over all $\s$ of $\sg\s$ times the product
\pagebreak

\[\prod_{i<j<m}(p+q\x_{\s(i)}\x_{\s(j)}-\x_{\s(i)})\;{1\ov (\xs{1}-1)\,(\xs{1}\xs{2}-1)\cd(\xs{1}\xs{2}\cd\xs{m-1}-1)}\]
\[\times\prod_{m\le i<j}(p+q\x_{\s(i)}\x_{\s(j)}-\x_{\s(i)})\;{\x_{\s(m+2)}\x_{\s(m+3)}^2\cd\x_{\s(N)}^{N-m}\ov (1-\x_{\s(m+1)}\cd\x_{\s(N)})\cd(1-\x_{\s(N)})}\]
\[\times \prod_{{i<m\le j}}(p+q\x_{\s(i)}\x_{\s(j)}-\x_{\s(i)}).\]

Now we take a fixed $U\subset\{1,\ldots,N\}$ with $|U|=m-1$, and sum over all permutations $\s$ such that $\s(i)\in U$ when $i\le m-1$. Notice that the last product above is equal to
\be\prod_{{i\in U\atop j\in U^c}}(p+q\x_i\x_j-\x_i),\label{last}\ee
and so is independent of the particular permutation. Therefore we may sum independently over bijective maps $\{1,\ldots,m-1\}\to U$ for the first factors and bijective maps $\{m,\ldots,N\}\to U^c$ for the second factors. 

To keep track of the signs of the permutations we observe that $\sg\s$ equals the product of the signs of the two restrictions of $\s$ times $\sg U$. So we eventually have to multiply by $\sg U$.
 
For the second factor above we use (\ref{identity1}) and find that the sum equals
\be p^{(N-m)(N-m+1)/2}\,\(1-\ds{\prod_{i\in U^c}\x_i}\)\,{\ds{\prod_{{i<j\atop i,j\in U^c}}(\x_j-\x_i)}\ov\ds{\prod_{i\in U^c}(1-\x_i)}}.\label{secondfactor}\ee

For the first factor we use (\ref{identity1}) and find that the sum equals
\be q^{(m-1)(m-2)/2}\; {\ds{\prod_{{i<j\atop i,j\in U}}(\x_j-\x_i)}\ov\ds{\prod_{i\in U}(\x_i-1)}}.\label{firstfactor}\ee

This is to be multiplied by (\ref{last}), divided by $\prod_{i<j}(p+q\x_i\x_j-\x_i)$, and then integrated over large contours $\C{R}$ with respect to the $\x_i$ with $i\in U$. But we want all integrals to be taken over the same contours $\C{r}$ so we want to replace the integral with all contours $\C{R}$ with a sum of integrals with all contours $\C{r}$.

In our application of Lemma 5.1 $N$ will be replaced by $m-1$,\ $\{1,\ld,N\}$ will replaced by $U$, and $S$ will be replaced by $T$, which is the reason we defined $\s(T,U)$ as the sum of the positions of the elements of $T$ in $U$. The coefficients become in this notation 

\be{q^{\s(\TS,\,U)-(m-1)\,|\TS|}\ov p^{\,|T|(m+|\TS|)/2-\s(T,U)}}.\label{coeff}\ee

If we put all integrands together the result is, aside from the powers of $p$ and $q$ in (\ref{firstfactor}) and (\ref{secondfactor}),
\be{\Big(1-\ds{\prod_i \x_i}\Big)\ds{\prod_{i<j}(\x_j-\x_i)}\ov\ds{\prod_{i<j}(p+q\x_i\x_j-\x_i)\,\prod_i(1-\x_i)}}\,\prod_i \(\x_i^{x-y_i-1}\,e^{\ep(\x_i)\,t}\)\quad(\textrm{all indices in}\ U^c)\label{Tc}\ee
\be\times {\ds{\prod_{i<j}(\x_j-\x_i)}\ov\ds{\prod_{i<j}(p+q\x_i\x_j-\x_i)\,\prod_i(1-\x_i)}}\,\prod_i \(\x_i^{x-y_i-1}\,e^{\ep(\x_i)\,t}\)\quad(\textrm{all indices in}\ U)\label{T}\ee
\[\times \ds{\prod_{{i>j\atop i\in U,\,j\in U^c}}{p+q\x_i\x_j-\x_i\ov
p+q\x_i\x_j-\x_j}}.\]
We apply Lemma~5.1 with  
\[g(\x)=\prod_{i\in U} \(\x_i^{x-y_i-1}\,e^{\ep(\x_i)\,t}\)
\;\int_{\C{r}}\cd\int_{\C{r}}\ds{\prod_{{i>j\atop i\in U,\,j\in U^c}}{p+q\x_i\x_j-\x_i\ov p+q\x_i\x_j-\x_j}}\,\prod_{j\in U^c}d\x_j,\]
where $g(\x)=g(\{\x_i\}_{i\in U})$.

The poles of the integrand defining $g$ are at
\[\x_i=(\x_j-p)/q\x_j,\]
where $j\in U^c$. Since $\x_j$ can be arbitrarily small the pole is outside $\C{R}$ for any $R$, so $g$ is analytic for all $\x_i \ne0$.

To check the main hypothesis on $g$ we observe that after the substitution $\x_i\to p/(1-q\x_k)$ the product $\x_i^{x-y_i-1}\,\x_k^{x-y_k-1}$ becomes $O(\x_k^{y_i-y_k})$ as $\x_k\to\iy$ which is $O(\x_k\inv)$ since $i<k$. The sum of the exponents in the last factor which involve $\x_i$ and $\x_k$ is
\[1-q\x_k+{pq\ov 1-q\x_k} +{p\ov \x_k}+q\x_k,\]
which is bounded at infinity. As for the integral, we show that it is bounded when one $\x_k\to\iy$ while the others are bounded. If we take a fixed $j>k$ with $j\in U^c$, the pole at $\x_j=p/(1-q\x_k)$ passes across the $\x_j$-contour when $\x_k\to\iy$. The residue of the factor in the product that contributes the pole is seen to be $O(1)$ as are the other factors involving $i\ne k$. We do this with each $j$ in turn and end up with a sum of integrals each of which is $O(1)$. If we make the substitution $\x_i\to p/(1-q\x_k)$ then $\x_k\to\iy$ while $p/(1-q\x_k)$ remains bounded, so $g$ satisfies the required hypothesis.

To find the summand in (\ref{intsum2}) we must evaluate $g(\x)$ where all the $\x_i$ with $i\in \TS$ set equal to 1. Each factor in the product in the integrand with such an $i$ is $(-q/p)$. For each such $i$ the number of $j$ in the product in the integrand satisfying $i>j$  and $j\in U^c$ equals $i$ minus the position of $i$ in $U$. The sum of this over all $i\in\TS$ equals $\s(\TS)-\s(\TS,\,U)$. If we multiply (\ref{coeff}) by $(-q/p)$ to this power the result may be written
\[ (-1)^{\s(\TS)-\s(U\bs T,U)}\,{q^{\s(\TS)-(m-1)\,|\TS|}\ov p^{\s(\TS)+
|T|(m+|\TS|)/2-m\,(m-1)/2}},\]
since $\s(U\bs T,U)+\s(T,U)=m(m-1)/2$.

There are also the factors $(-1)^{|U\bs T|}$ coming from (\ref{intsum2}) and $(-1)^{|U|}$ coming from the fact that $\prod(\x_i-1)$ appears in (\ref{firstfactor}) rather than $\prod(1-\x_i)$. These factor combine as $(-1)^{|T|}$.

For the integrand we must combine (\ref{Tc}), (\ref{T}) with the $\x_i$ with $i\in\TS$ set equal to 1, and the integrand in $g$ with the $\x_i$ with $i\in\TS$ set equal to 1. The result is (\ref{YTU}) with $U$ replaced by $U^c$. If we take account of the original factors in (\ref{firstfactor}) and (\ref{secondfactor}) we have established  Lemma 5.2.

\noi{\bf Proof of Theorem 5.1}. Suppose temporarily that $q\ne0$ also. We take a fixed $S\subset\{1,\ld,N\}$ with $|S^c|\le m$ and in Lemma~5.2 sum over all $T$ and $U$ with $T\subset U$ and $T\cup U^c=S$. Let us write everything in terms of $T$ and $U^c$.

For any $i\in U$, the position of $i$ in $U$ equals $\#\{j:j\le i,\ j\in U\}$, so
\[\s(T,U)=\#\{(i,j):i\ge j,\ i\in T,\,j\in U\}.\]
In particular,
\[\s(\TS,\,U)=\#\{(i,j):i\ge j,\ i\in \TS,\,j\in U\},\]
\[\s(\TS)=\#\{(i,j):i\ge j,\ i\in \TS,\,j\in \{1,\cd,N\}\},\]
and so
\[\s(\TS)-\s(\TS,\,U)=\#\{(i,j):i\ge j,\ i\in \TS,\,j\in U^c\}.\]
Also, 
\[\sg U=(-1)^{\#\{(i,j)\,:\,i>j,\;i\in U,\,j\in U^c\}}.\]
Thus 
\be(-1)^{|T|+\#\{(i,j)\,:\,i>j,\;i\in T,\,j\in U^c\}}.\label{combined}\ee
is the combined power of $-1$ that occurs.

The powers of $p$ and $q$ that occur in the summation are, since $U\bs T=S^c$,
\be{q^{\s(S^c)-(m-1)\,|S^c|}\ov p^{\s(S^c)+|T|(m+|S^c|)/2}},\label{pqpower}\ee
and only depends on $|T|$, given $S$. 

In (\ref{YTU}), we write
\[\ds{\prod_{{i<j\atop i,j\in U^c\ {\rm or}\ i,j\in T}}(\x_j-\x_i)}
={\ds{\prod_{{i<j\atop i,j\in U^c\cup T}}(\x_j-\x_i)}\ov 
\ds{\prod_{{i<j\atop i\in U^c,\;j\in T}}(\x_j-\x_i)}\,\ds{\prod_{{i<j\atop i\in T,\;j\in U^c}}(\x_j-\x_i)}}.\]
The denominator may be written
\[\ds{\prod_{{i>j\atop i\in T,\;j\in U^c}}(\x_i-\x_j)}\,\ds{\prod_{{i<j\atop i\in T,\;j\in U^c}}(\x_j-\x_i)}=
(-1)^{\#\{(i,j):i>j,\;i\in T,\,j\in U^c\}}\prod_{i\in T,\;j\in U^c}(\x_j-\x_i).\]
The power of $-1$ combined with (\ref{combined}) equals $(-1)^{|T|}$ and therefore our integrand, aside from this power of $-1$, may be written
\[\Big(1-\prod_{i\in U^c} \x_i\Big){\ds{\prod_{i\in T,\;j\in U^c}(p+q\x_i\x_j-\x_i)}
\ov \ds{\prod_{i\in T,\,j\in U^c}(\x_j-\x_i)}}\,{\ds{\prod_{i<j}(\x_j-\x_i)}\ov\ds{\prod_{i}(1-\x_i)\,\prod_{i<j}(p+q\x_i\x_j-\x_i)}}
\prod_i \Big(\x_i^{x-y_i-1}\,e^{\ep(\x_i)\,t}\Big),\]
where indices not specified range over $T\cup U^c$.

Now we take a fixed $S\subset\{1,\ld,N\}$ with $|S^c|<m$ and in Lemma~5.2 first sum over all $T$ and $U$ such that $T\cup U^c=S$. The condition $|U|=m-1$ translates to $|T|=m-1-|S^c|$. Apply (\ref{identity2}) with $\{1,\ld,N\}$ replaced by $S$, with $S$ replaced $T$, and with $m$ replaced by $m-1$. We obtain (\ref{Th5.1}), and this completes the proof when $q\ne0$. 

We can remove this condition by taking the $q\to0$ limit since no pole tends to zero when $q\to0$.

We now obtain the expansion where all integrals are taken over large contours.
\pagebreak

\noi{\bf Theorem 5.2}. We have when $q\ne0$
\[\P(x_{m}(t)=x)=(-1)^{m+1}(pq)^{m(m-1)/2}\]
\be\times\sum_{|S|\ge m}\br{|S|-1}{|S|-m}\,{p^{\s(S)- m\,|S|}\ov q^{
\s(S)-|S|(|S|+1)/2}}\int_{\C{R}}\cd\int_{\C{R}}I(x, Y_S,\x)\,d^{|S|}\x,\label{Th5.2}\ee
where $R$ is so large that the poles of the integrand lie inside $\C{R}$. The sum is taken over all subsets $S$ of $\{1,\ld,N\}$ with $|S|\ge m$.

\noi{\bf Proof}.  Denote by $\tl\P$ the probabilities for the process with $p$ and $q$ interchanged. As in the remark following the proof of Theorem~2.1, $\P(x_m(t)=x)$ is equal to $\tl\P(x_{N-m+1}(t)=-x)$  with initial configuration $Y$ replaced by $\{-y_N,\ld,-y_1\}$. In the integrals in (\ref{Th5.2}) we make the replacements $\x_i\to 1/\x_{N-i+1}$. The upshot is that in (\ref{Th5.1}) we replace $m$ by $N-m+1$, in the coefficients we replace $S$ by $\tl S$, we multipy by $(-1)^{|S|+1}$ (because of the sign change in the integrand), and take the integrals over $\C{R}$. A little algebra, using the general fact
\[\s(\tl S)=\sum_{i\in S}(N-i+1)=(N+1)\,|S|-\s(S),\]
shows that the result is (\ref{Th5.2}). This complete the proof of Theorem~5.2.

As with the first particle when $Y$ is infinite and bounded below, we can show that the sum (\ref{Th5.2}) converges when it is taken over all finite subsets of $\Z^+$. This gives the probability for infinitely many particles. 

In the special case $Y=\Z^+$ we can evaluate the sum over all sets of a given cardinality. We define the $k$-dimensional integrand
\[J_k(x,\x)=\prod_{i\ne j}{\x_j-\x_i\ov p+q\x_i\x_j-\x_i}\;{1-\x_1\cd\x_k\ov\ds{\prod_i(1-\x_i)\,(q\x_i-p)}}\,\prod_i\(\x_i^{x-1}e^{\ep(\x_i)t}\),\]
where all indices run over $\{1,\ld,k\}$. The result we obtain is

\noi{\bf Corollary}. When $Y=\Z^+$ we have when $q\ne0$
\[\P(x_m(t)=x)=(-1)^{m+1}\,q^{m(m-1)/2}\]
\be\times \sum_{k\ge m}{1\ov k!}\,\br{k-1}{k-m}\,p^{(k-m)(k-m+1)/2}\;q^{k(k+1)/2}\int_{\C{R}}\cd\int_{\C{R}}J_k(x,\x)\,d\x_1\cd d\x_k.\label{Cor}\ee

\noi{\bf Proof}. We sum first over all $S\subset\Z^+$ with $|S|=k$. If $S=\{z_1,\ld,z_k\}$ we make the variable changes $\x_{z_1}\to\x_1,\ld,\x_{z_k}\to\x_k$ in the integral of $I(x, Y_S,\x)$, so all integration are over the variables $\x_i$ with $i\in\{1,\ld,k\}$. The integrand becomes 
\[\prod_{i<j}{\x_j-\x_i\ov p+q\x_i\x_j-\x_i}\;{1-\x_1\cd\x_k\ov\ds{\prod_i(1-\x_i)}}\,\prod_i\(\x_i^{x-z_i-1}\,e^{\ep(\x_i)t}\),\]
where all indices run over $\{1,\ld,k\}$. The only part of the coefficient in (\ref{Th5.2}) involving more than $|S|$ is $(p/q)^{\s(S)}$. When we multiply by this the product may be written
\[{\ds{\prod_{i<j}(\x_j-\x_i)}\ov \ds{\prod_{i\ne j}(p+q\x_i\x_j-\x_i)}}\;{1-\x_1\cd\x_k\ov\ds{\prod_i(1-\x_i)}}\,\prod_i\(\x_i^{x-1}\,e^{\ep(\x_i)t}\)\]
\[\times \prod_{i>j}(p+q\x_i\x_j-\x_i)\;\prod_i\Big({q\ov p}\,\x_i\Big)^{-z_i}.\]
If we sum over all $\{z_i\}$ with $0<z_1<\cd<z_k$ the last product becomes
\[\prod_{i>j}(p+q\x_i\x_j-\x_i)\;{1\ov \Big(({q\ov p}\,\x_1)
({q\ov p}\,\x_2)\cd({q\ov p}\,\x_k)-1\Big)\,\Big
(({q\ov p}\,\x_2)\cd({q\ov p}\,\x_k)-1\Big)\cd\Big(({q\ov p}\,\x_k)-1\Big)}.\]

The integral is unchanged if we antisymmetrize this, because all other factors are symmetric except for the Vandermonde. We can bring this antisymmetrization to the form of identity (\ref{identity1'}) if we make the substitutions
\[\x_i={p\ov q}\,\e_{k-i+1}.\]
The second factor becomes
\[{1\ov (\e_1-1)\,(\e_1\,\e_2-1)\cd (\e_1\,\e_2\cd \e_k-1)},\]
while the first factor (after the index changes $i\to k-i+1,\ j\to k-j+1$) becomes
\[\({p\ov q}\)^{k(k-1)/2}\,\prod_{i<j}(q+p\e_i\e_j-\e_i).\]
Now we can apply (\ref{identity1'}) with $p$ and $q$ interchanged and we see that the antisymmetrization of this is
\[{1\ov k!}\,{p^{k(k-1)}\ov q^{k(k-1)/2}}\,{\ds{\prod_{i<j}(\e_j-\e_i)}\ov\ds{\prod_i(\e_i-1)}}
={1\ov k!}\,p^{k(k+1)/2}\,{\ds{\prod_{i>j}(\x_j-\x_i)}\ov\ds{\prod_i(q\x_i-p)}}.\]
If we recall the remaining factor
\[{p^{- m\,|S|}\ov q^{-|S|(|S|+1)/2}}={q^{k(k+1)/2}\ov p^{mk}}\]
in the coefficient in (\ref{Th5.2}) we see that we obtain formula (\ref{Cor}).

\noi{\bf Remark}. The power of $p$ on the right side of (\ref{Cor}) is always nonnegative and is zero only when $k=m$. Hence when $p=0$, in other words in the TASEP where particles move to the left, only one term survives and we obtain
\[\P(x_{m}(t)=x)={(-1)^{m(m-1)/2}\ov m!}
\int_{\C{R}}\cd\int_{\C{R}}\prod_{i<j}(\x_j-\x_i)^2\,
{\x_1\cd\x_m-1\ov \ds{\prod_i(\x_i-1)^m}}\]
\[\times\prod_i\Big(\x_i^{x-m-1}\,e^{(\x_i-1) t}\Big)\,d\x_1\cd d\x_m.\]
For $\P(x_{m}(t)\le x)$ we sum $\P(x_{m}(t)=y)$ over all $y\le x$ (which we may, since $R>1$), obtaining
\[\P(x_{m}(t)\le x)={(-1)^{m(m-1)/2}\ov m!}
\int_{\C{R}}\cd\int_{\C{R}}\prod_{i<j}(\x_j-\x_i)^2\;
\prod_i(\x_i-1)^{-m}\]
\[\times\prod_i\Big(\x_i^{x-m}\,e^{(\x_i-1) t}\Big)\,d\x_1\cd d\x_m.\]
By a general identity \cite{A} this equals
\[(-1)^{m(m-1)/2}\,\det\(\int_{\C{R}}\x^{i+j+x-m}\,(\x-1)^{-m}\,e^{(\x-1)t}\,d\x\)_{i,\,j=0,\ld,m-1}.\]
After reversing the order of the columns this becomes a Toeplitz determinant equal to the determinant of 
R\'akos-Sch\"utz \cite[(12)]{RS} which they used to obtain Johansson's result \cite{Jo}.

\setcounter{equation}{0}\renewcommand{\theequation}{6.\arabic{equation}}
 
\begin{center}{\bf VI. Proofs of the Identities}\end{center}

\noi{\bf Proof of identity (\ref{identity1})}.\footnote{When we showed Doron Zeilberger the identity when it was still a conjecture he suggested \cite{Z} that problem VII.47 of \cite{PSz}, an identity of I.~Schur, had a similar look about it and might be proved in a similar way. He was right.} We use induction on $N$, and assume that the identity holds for $N-1$. (It clearly holds for $N=1$.) Call the left side $\ph_N(\x_1,\ldots,\x_N)$. We first sum over all permutations such that $\s(1)=k$, and then sum over $k$. If we observe that the inequality $i<j$ becomes $j\ne i$ when $i=1$, we see that what we get for the left side is
\[{1\ov1-\x_1\,\x_2\cd \x_N}\sum_{k=1}^N(-1)^{k+1}\prod_{j\ne k}(p+q\x_k\x_j-\x_k)\cdot\prod_{j\ne k}\x_j\cdot\ph_{N-1}(\x_1,\ldots,\x_{k-1},\x_{k+1},\ldots,\x_N),\]
which may also be written
\[{\x_1\,\x_2\cd \x_N\ov1-\x_1\,\x_2\cd \x_N}\sum_{k=1}^N(-1)^{k+1}\,\prod_{j\ne k}(p+q\x_k\x_j-\x_k)\cdot\x_k\inv\,\ph_{N-1}(\x_1,\ldots,\x_{k-1},\x_{k+1},\ldots,\x_N).\]
We want to show that this equals the right side of (\ref{identity1}), and the induction hypothesis gives
\[\ph_{N-1}(\x_1,\ldots,\x_{k-1},\x_{k+1},\ldots\x_N)=p^{(N-1)(N-2)/2} \;{\ds{\prod_{i<j;\,i,j\ne k}}(\x_j-\x_i)\ov \ds{\prod_{j\ne k} (1-\x_j)}}.\]
After some multiplying, dividing, and computing powers of $-1$ we see that what we want to show is 
\be\sum_{k=1}^N\prod_{j\ne k}(p+q\x_k\x_j-\x_k)\cdot{1-\x_k\ov\x_k}{1\ov\prod_{j\ne k}(\x_j-\x_k)}=p^{N-1}{1-\x_1\,\x_2\cd \x_N\ov\x_1\,\x_2\cd \x_N}.\label{2}\ee
 
If we change the first product on the left to run over all $j$ we have to divide by $p+q\x_k^2-\x_k$. Setting $p=1-q$ shows that
\[{1-\x_k\ov p+q\x_k^2-\x_k}={1\ov p-q\x_k}.\]
So the left side of (\ref{2}) equals
\[\sum_{k=1}^N\,\prod_{j=1}^N(p+q\x_k\x_j-\x_k)\cdot{1\ov
\x_k\,(p-q\x_k)}\,{1\ov\prod_{j\ne k}(\x_j-\x_k)}.\]

We evaluate this by integrating  
\[\prod_{j=1}^N(p+qz\x_j-z)\cdot{1\ov
z\,(p-qz)}\cdot{1\ov\prod_{j=1}^N(\x_j-z)}\]
over a large circle. Since the integrand is $O(z^{-2})$ for large $z$ the integral is zero. There are poles at 0 and the $\x_k$, and the sum of the residues there is
\[{p^{N-1}\ov \prod_j\x_j}-\sum_{k=1}^N\prod_{j=1}^N(p+q\x_k\x_j-\x_k)\cdot{1\ov
\x_k\,(p-q\x_k)}\,{1\ov \prod_{j\ne k}(\x_j-\x_k)}.\]
There is also a pole at $z=p/q$ and for the residue there we compute
\[{p+p\x_j-p/q\ov\x_j-p/q}=p\,{q+q\x_j-1\ov q\x_j-p}=p,\]
so the residue at $p/q$ is $-p^{N-1}$. This gives
\[\sum_{k=1}^N\,\prod_{j=1}^N(p+q\x_k\x_j-\x_k)\cdot{1\ov
\x_k\,(p-q\x_k)}\,{1\ov \prod_{j\ne k}(\x_j-\x_k)}={p^{N-1}\ov \prod_j\x_j}-p^{N-1},\]
as desired.

This completes the proof of identity (\ref{identity1}).

\noi{\bf Proof of identity (\ref{identity2})}. We use the easily established recursion formula
\be \br{N}{m}=p^m\,\br{N-1}{m}+q^{N-m}\,\br{N-1}{m-1}\label{recursion}\ee
and a preliminary simpler identity,
\be\sum_{|S|=m} \prod_{{ i\in S\atop  j\in S^c}} {p+q\x_i\x_j-\x_i\ov \x_j-\x_i}=\br{N}{m},\label{identity3}\ee
where $S$ is as before.

We prove this\footnote{The proof is a modification of one found by Anne Schilling \cite{AS} for an equivalent identity.} by induction on $N$, so assume (\ref{identity3}) holds for $N-1$. The left side is symmetric in the $\x_i$ and is $O(1)$ as any $\x_i\to\iy$ with the other $\x_j$ fixed. If we multiply the left side by the Vandermonde $\prod_{i<j}(\x_i-\x_j)$ we obtain an antisymmetric polynomial which is $O(\x_i^{N-1})$ as any $\x_i\to\iy$ with the others fixed, so it has degree at most $N-1$ in each $\x_i$. Being antisymmetric it is divisible by the Vandermonde, and having degree at most $N-1$ in each of the $\x_i$ separately it must be a constant times the Vandermonde. Thus the left side of (\ref{identity3}) is a constant, say $C_{N,m}$.

To evaluate the constant set $\x_N=1$. For convenience we write
\[{p+q\x_i\x_j-\x_i\ov \x_j-\x_j}=U(\x_i,\x_j).\]
We have
\[C_{N,m}=\sum_{S} \prod_{{ i\in S\atop  j\in S^c}} U(\x_i,\x_j)\Big|_{\x_1=1}=
q^{N-m}\,\sum_{N\in S}\prod_{{i\in S\bs N\atop j\in S^c}}U(\x_i,\x_j)
+p^m\,\sum_{N\not\in S}\prod_{{i\in S\atop j\in S^c\bs N}}U(\x_i,\x_j).\]
By the induction hypothesis the right side equals 
\[q^{N-m}\,\br{N-1}{m-1}+p^m\,\br{N-1}{m},\]
and this equals $\br{N}{m}$ by (\ref{recursion}).
This establishes identity (\ref{identity3}).

The proof of (\ref{identity2}) runs along the same lines. We interpret both sides to be zero when $N=m$, and do an induction on $N\ge m$. So we assume $N>m$ and that the formula holds for $N-1\ge m-1$. We quickly deduce that the left side is a polynomial of degree at most one in each $\x_i$. If we call the left side $C_{N,m}(\x)$ then
\[C_{N,m}(\x_1,\ld,\x_{N-1},1)=q^{N-m}\,\sum_{N\in S}\prod_{{i\in S\bs N\atop j\in S^c}}U(\x_i,\x_j)\cdot\Big(1-\prod_{j\in S^c}\x_j\Big)\]
\[+p^m\,\sum_{N\not\in S}\prod_{{i\in S\atop j\in S^c\bs N}}U(\x_i,\x_j)\cdot\Big(1-\prod_{j\in S^c\bs N}\x_j\Big).\]
\[=q^{N-m}\,C_{N-1,m-1}(\x_1,\ld,\x_{N-1})+p^m\,C_{N-1,m}(\x_1,\ld,\x_{N-1}).\]
Similar relations hold for the other $\x_i$. Notice that when $N=m$ the second sum above does not appear and $C_{m-1,m-1}(\x)=0$ so this is consistent with our initial condition.

If we call the right side of (\ref{identity2}) $C_{N,m}'(\x)$ we see from (\ref{recursion}) that the same relations hold for $C_{N,m}'(\x)$ (with initial condition $C_{m,m}'(\x)=0$) and so for the difference $D_{N,m}(\x)=C_{N,m}(\x)-C_{N,m}'(\x)$. The induction hypothesis gives $D_{N-1,m}=D_{N-1,m-1}=
\nolinebreak 0$, so we have shown that for any $i$
\[D_{N,m}(\x)|_{\x_i=1}=0.\]
Any polynomial which has degree at most one in each $\x_i$ and vanishes when any $\x_i=1$ is a constant times $\prod (\x_i-1)$,\footnote{Such a polynomial must be of the form $\x_N-1$ times a polynomial in $\x_1,\ld,\x_{N-1}$ with the same property, so the statement follows by induction.} so $D_{N,m}(\x)$ has this form.

We have shown that
\[C_{N,m}(\x)=C_{N,m}'(\x)+c\,\prod_i(\x_i-1)\]
for some $c$. We show that $c=0$ by computing asymptotics as $\x_N\to\iy$. All terms are asymptotically a constant times $\x_N$. If in the sum in (\ref{identity2}) $N\in S$ then the corresponding summand is $O(1)$. So we need consider only those $S$ for which $N\not\in S$. In the product in the summand, if $j=N$ then the corresponding product over $i$ has the limit $q^m\,\prod_{i\in S}\x_i$ since there are $m$ factors with limit $q\x_i$. It follows that 
\[\lim_{\x_N\to\iy}{C_{N,m}(\x)\ov \x_N}=-q^m\,\sum_{{|S|=m\atop S\subset \{1,\ld,N-1\}}} \prod_{{ i\in S\atop  j\in S^c, j<N}} {p+q\x_i\x_j-\x_i\ov \x_j-
\x_i}\cdot\prod_{i\in S}\x_i\cdot \prod_{j<N,\,j\in S^c}\x_i\]
\[=-q^m\,\prod_{i<N}\x_i\;\sum_{{|S|=m\atop S\subset \{1,\ld,N-1\}}} \prod_{{ i\in S\atop  j\in S^c, j<N}} {p+q\x_i\x_j-\x_i\ov \x_j-
\x_i}.\]
Identity (\ref{identity3}) tells us that this equals
\[-q^m\,\prod_{i<N}\;\x_i\br{N-1}{m}.\]
Clearly $C_{N,m}'(\x)$ has the same asymptotics, so $c=0$.

\begin{center}{\bf Acknowledgments}\end{center}

The authors thank Thomas Liggett, Anne Schilling and Doron Zeilberger for their invaluable assistance. Helpful communications with Alex Gamburd and Gunter Sch\"utz are gratefully acknowledged as are comments by Timo Sepp\"al\"ainen and Pavel Bleher.  This work was supported by the National Science Foundation under grants DMS--0553379 (first author) and DMS--0552388 (second author).


\begin{thebibliography}{99}


\bibitem{ADHR} Alcaraz, F.C,  Droz, M.,  Henkel, M., and Rittenberg, V.:
Reaction-diffusion processes, critical dynamics, and quantum chains.
Ann.\  Physics \textbf{230}, 250--302  (1994).

\bibitem{A} Andr{\'e}ief, C.: Note sur une relation
les int\'egrales d\'efinies des produits des fonctions,
M\'em.\ de la Soc.\ Sci., Bordeaux {\bf 2}, 1--14 (1883).

\bibitem{BDJ} Baik, J., Deift, P.A., Johansson, K.: On the distribution of
the length of the longest increasing subsequence in a random permutation.
J.\ Am.\ Math.\ Soc.\ \textbf{12}, 1119--1178 (1999).

\bibitem{GS} Gwa, L.-H., Spohn, H.:  Bethe solution for the dynamical-scaling
exponent of the noisy Burgers equation.  Phys.\ Rev.\ A \textbf{46}, 844--854  (1992).

\bibitem{HKPV} Hough, J.B., Krishnapur, M., Peres, Y., Vir\'ag, B.: Determinantal processes
and independence. Probability Surveys \textbf{3}, 206--229 (2006).

\bibitem{Jo} Johansson, K.: Shape fluctuations and random matrices.
Commun.\ Math.\ Phys. \textbf{209}, 437--476  (2000).

\bibitem{LL} Lieb, E.H., Liniger, W.: Exact analysis of an interacting Bose gas. I.
 The general solution and the ground state. Phys.\ Rev.\ \textbf{130}, 1605--1616 (1963).

\bibitem{Li1} Liggett, T.M.: \textit{Interacting Particle Systems}.
[Reprint of the 1985 original.]  Berlin, Springer-Verlag,
 2005.

\bibitem{Li2} Liggett, T.M.: \textit{Stochastic Interacting Systems: Contact, Voter and
Exclusion Processes}. Berlin, Springer-Verlag, 1999.

\bibitem{Li3} Liggett, T.M.: Private communication April 17, 2007.

\bibitem{PSz} P\'olya, G., Szeg\"o, G.: \textit{Aufgaben und Lehrs\"atze aus der Analysis}. Berlin, Springer-Verlag, 1964.

\bibitem{PS} Pr\"ahofer, M., Spohn, H.: Current fluctuations for the totally
asymmetric simple exclusion process. \textit{In and Out of Equilibrium},
Progress in Probability \textbf{51}, 185--204  (2000).

\bibitem{RS} R\'akos, A., Sch\"utz, G.M.: Current distribution and random matrix
ensembles for an integrable asymmetric fragmentation process.
J.\ Stat.\ Physics \textbf{118}, 511--530  (2005).

\bibitem{Sa} Sasamoto, T.: Spatial correlations of the 1D KPZ surface on a flat substrate. J. Phys. A \textbf{38}, L549--L556 (2005).

\bibitem{AS} Schilling, A.: Private communication Apr. 6, 2007.

\bibitem{Sc}  Sch\"utz, G.M.: Exact solution of the master equation for the
asymmetric exclusion process. 
J.\ Stat.\ Physics \textbf{88}, 427--445 (1997). 

\bibitem{Sos} Soshnikov, A.: Determinantal random fields. Russian Math.\ Surveys
\textbf{55}, 923--975 (2000).

\bibitem{Spo} Spohn, H.: Exact solutions for KPZ-type growth processes,
random matrices, and equilibrium shapes of crystals. Physica A \textbf{369}, 71--99 (2006).

\bibitem{Spi} Spitzer, F.:  Interaction of Markov processes.  Adv.\ Math.\
\textbf{5}, 246--290  (1970).

\bibitem{Sut} Sutherland, B.: \textit{Beautiful Models: 70 Years of Exactly Solvable
Quantum Many-Body Problems}, World Scientific, 2004.

\bibitem{Yang2}  Yang, C.N., Yang, C.P.:  One-dimensional chain of
anisotropic spin-spin interactions. I. Proof of Bethe's hypothesis for the ground state
in a finite system.  Phys.\ Rev.\ \textbf{150},  321--327  (1966).

\bibitem{Yau} Yau, H.-T.:  $(\log t)^{2/3}$ law of the two dimensional asymmetric simple
exclusion process. Ann.\ Math.\ \textbf{159}, 377--405 (2004).

\bibitem{Z} Zeilberger, D.: Private communication Feb. 14, 2007.

\end{thebibliography}
\end{document}